\documentclass[final,12pt,3p]{elsarticle}

\usepackage{graphicx}
\usepackage{amssymb}
\usepackage{amsthm}
\usepackage{amsmath}
\usepackage{algorithm}
\usepackage{algorithmic}
\usepackage{tikz}
\usepackage{lineno}
\usepackage{verbatim} 
\usepackage{hyperref}
 \usepackage{xcolor}
\usepackage{caption}
\usetikzlibrary{shapes}
\usepackage{subcaption}
\usepackage{geometry}
\usepackage{mathrsfs} 
\usepackage{array}
\usepackage{multirow}
\usepackage{booktabs} 
\usepackage{adjustbox}
\usepackage[title]{appendix}

\biboptions{sort&compress}
\hypersetup{
colorlinks=true,
linkcolor=red,
filecolor=green,
urlcolor=green,
anchorcolor=green,
citecolor=cyan,
pdftitle={Overleaf Example},
pdfpagemode=FullScreen,
}
\begin{document}
	
\begin{tikzpicture}[remember picture,overlay]
	\node[anchor=north east,inner sep=20pt] at (current page.north east)
	{\includegraphics[scale=0.2]{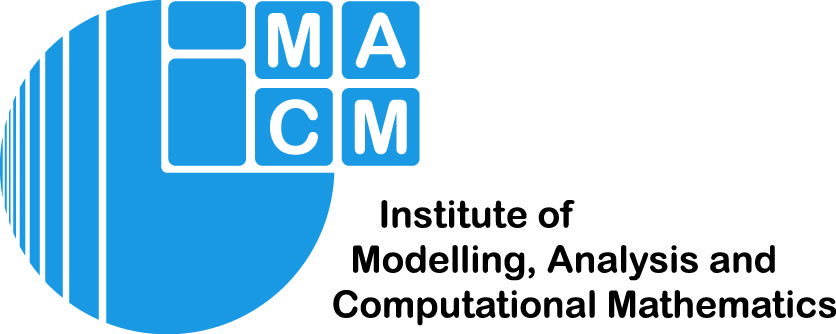}};
\end{tikzpicture}

\begin{frontmatter}

\title{Constraint-Aware Physics-Informed Neural Networks for SEIR Reaction-Diffusion Epidemic Models with Vital Dynamics}

\author[BUW,MAIS]{Achraf Zinihi}
\ead{a.zinihi@edu.umi.ac.ma} 

\author[BUW]{Matthias Ehrhardt\corref{Corr}}
\cortext[Corr]{Corresponding author}
\ead{ehrhardt@uni-wuppertal.de}

\address[BUW]{University of Wuppertal, Applied and Computational Mathematics,\\
Gaußstrasse 20, 42119 Wuppertal, Germany}

\address[MAIS]{Department of Mathematics, AMNEA Group, Faculty of Sciences and Techniques,\\
Moulay Ismail University of Meknes, Errachidia 52000, Morocco}


\begin{abstract}
Reaction-diffusion epidemic models with vital dynamics are an important framework for describing the spatial and temporal spread of infectious diseases. 
In this work, we present a constraint-aware, physics-informed neural network (PINN) approach to an SEIR reaction-diffusion system with homogeneous Neumann boundary conditions. 
Due to the scarcity of spatial epidemiological datasets, we generate synthetic benchmark data using structure-preserving implicit-explicit nonstandard finite difference (NSFD) schemes that ensure positivity, boundedness, and numerical stability.

The PINN framework integrates PDE residuals, observational data, boundary conditions, and epidemiological constraints within a unified optimization procedure. 
Specifically, the loss function incorporates the non-negativity of compartment populations and the admissibility of epidemiological parameters.
We apply the method to forward simulation and inverse parameter estimation in one- and two-dimensional settings.
Numerical experiments demonstrate the framework's ability to accurately reconstruct spatiotemporal epidemic dynamics and reliably identify parameters, even when data is sparse or noisy. 
These results underscore the potential of constraint-aware PINNs as a robust, data-driven methodology for spatial epidemic modeling.

\end{abstract}


\journal{} 



%

\end{frontmatter}


\section{Introduction}\label{S1}
The spatiotemporal spread of infectious diseases is a primary concern in both mathematical epidemiology and public health policy. This concern has motivated the development of models capable of predicting, explaining, and ultimately controlling outbreaks. 
Classical compartmental models, such as the SIR, SIQR, and SEIR systems, have provided essential insights into epidemic dynamics since the pioneering work of Kermack and McKendrick \cite{Kermack1927, Hethcote2000}. 
These models are traditionally formulated as systems of ordinary differential equations (ODEs) and capture average population-level dynamics, but they neglect spatial heterogeneity and human mobility. Both of these factors are critical drivers of epidemic spread. 
To address these limitations, reaction-diffusion formulations of epidemic models have been introduced to simulate spatially structured outbreaks, wave-like infection fronts, and localized clustering phenomena \cite{Zinihi2025FDE, Wang2012, Wang2021, Zinihi2025OC}.

Despite their conceptual appeal, reaction-diffusion epidemic partial differential equations (PDEs) present substantial challenges in practice. 
Classical numerical methods, such as finite difference schemes, finite element methods, and spectral solvers \cite{Chang2022, Strikwerda2004, Zinihi2025S}, provide accurate approximations, but they suffer from the curse of dimensionality and become prohibitively expensive in two or more spatial dimensions. 
Furthermore, inverse problems, such as estimating epidemiological parameters 
from noisy, sparse, or heterogeneous surveillance data, are notoriously ill-posed \cite{Banks1989}, often yielding unstable or biologically implausible estimates.

In recent years, physics-informed neural networks (PINNs) have emerged as a promising approach to data-driven partial differential equation (PDE) modeling \cite{Raissi2019, Karniadakis2021}. 
PINNs embed governing equations and boundary conditions directly into the loss function of a neural network, enabling mesh-free solution of PDEs and the simultaneous integration of data and mechanistic knowledge.
This approach has grown rapidly across disciplines, including fluid mechanics, materials science, and epidemiology \cite{Kissas2020, Cuomo2022, Toscano2025}.
In the context of epidemic systems, PINNs have been successfully applied to ODE-based models of diseases such as COVID-19 and other respiratory viruses \cite{Zinihi2025MC, Berkhahn2022, Millevoi2024}, as well as for forecasting hospitalizations \cite{Nelson2024, Shamsara2025} and even for fractional-order epidemic dynamics \cite{Zinihi2025PINN, Igiri2025}. 
Recent works have emphasized bi-objective optimization \cite{Heldmann2023}, 
physics-informed graph neural ODEs \cite{Cheng2024}, neural ODE–Bayesian hybrids \cite{Liu2024}, 
and variance-aware PINNs \cite{Shan2025}, all of which point toward the rapidly evolving frontier of physics-informed machine learning in epidemiology.

Nevertheless, the application of PINNs to spatial epidemic PDEs is still in its early stages. 
Most existing works focus on non-spatial ODEs or simplified 1D PDE settings \cite{Madden2024, Chen2025}.  
Furthermore, the systematic enforcement of epidemiological constraints, such as the non-negativity of populations, the boundedness of parameters, and the conservation of the total population under vital dynamics, has rarely been addressed despite being essential for biologically meaningful results. 
Without these constraints, PINN approximations may yield unrealistic solutions, such as negative compartment sizes or unbounded growth.

Based on the above discussion, several important gaps remain in the current literature on PINN-based epidemic modeling. Most existing studies are limited to ODE systems or 1D PDE formulations. 
Meanwhile, 2D reaction–diffusion models, which are better suited for realistic spatial epidemic propagation, have not been widely explored. 
Additionally, epidemiological constraints, such as the non-negativity of compartment populations, bounded transmission rates, and population balance under vital dynamics, are rarely explicitly incorporated into the learning framework despite their importance for biological realism. 

Another major limitation is the lack of high-quality spatiotemporal benchmark datasets for reaction-diffusion epidemic PDEs. 
This makes systematic validation and reproducibility difficult. 
Consequently, synthetic data generation using structure-preserving methods, such as NSFD schemes, is essential. 
Finally, rigorous comparative studies of 1D and 2D epidemic PINN models with respect to predictive accuracy, computational efficiency, and parameter identifiability are lacking.

This paper presents a constraint-aware PINN framework for SEIR reaction-diffusion epidemic models with vital dynamics and homogeneous Neumann boundary conditions. 
This approach incorporates biologically motivated constraints into the loss function to enforce the positivity and boundedness of compartment populations, the admissibility of epidemiological parameters, and population balance under births and natural deaths. 
Synthetic 1D and 2D benchmark datasets are generated using structure-preserving implicit-explicit nonstandard finite difference (NSFD) schemes under various noise and sparsity conditions to enable systematic and reproducible evaluation. 
Furthermore, we perform comprehensive comparisons between 1D and 2D PINN formulations to assess predictive accuracy, robustness, and computational efficiency. 

%

The remainder of this work is structured as follows. 
Section~\ref{S2} introduces the proposed SEIR reaction-diffusion model with vital dynamics and establishes its well-posedness. 
Section~\ref{S3} details the construction of synthetic datasets using NSFD schemes and their numerical analysis. 
Section~\ref{S4} presents the constraint-aware PINN framework, including the network architecture, loss function design, and training strategy. 
Section~\ref{S5} reports numerical results, including forward predictions, inverse parameter estimation, and comparisons between 1D and 2D settings. 
Finally, Section~\ref{S6} concludes the paper with a summary of the findings.

\section{Mathematical Model}\label{S2}
We consider an SEIR-type epidemic system with \emph{vital dynamics} and spatial diffusion, defined on a bounded spatial domain $\Omega \subset\mathbb{R}^d$ $(d=1,2)$ with an outward unit normal vector $\vec{n}$ on the boundary $\partial \Omega$. 
The temporal domain is $t\in[0,T]$ for some final time $T>0$. 
The state variables are the population densities of the susceptible $S(t,x)$, exposed $E(t,x)$, infected $I(t,x)$, and recovered $R(t,x)$ classes at time $t$ and spatial position $x \in \Omega$.

\subsection{The PDE System}\label{S2.1}

\begin{figure}[H]
\centering
\begin{tikzpicture}[node distance=3.5cm]
\node (S) [rectangle, draw, minimum size=1cm, fill=cyan!30] {S};
\node (E) [rectangle, draw, minimum size=1cm, fill=orange!40, right of=S] {E};
\node (I) [rectangle, draw, minimum size=1cm, fill=red!40, right of=E] {I};
\node (R) [rectangle, draw, minimum size=1cm, fill=green!30, right of=I] {R};
\draw [->, thick, >=latex, line width=1pt] (-1.5,0) -- ++(S) node[midway,above]{$\Lambda$};
\draw [->, thick, >=latex, line width=1pt] (S) -- (E) node[midway,above]{$\beta p S I$};
\draw [->, thick, >=latex, line width=1pt] (E) -- (I) node[midway,above]{$\delta E$};
\draw [->, thick, >=latex, line width=1pt] (I) -- (R) node[midway,above]{$\gamma I$};
\draw [->, thick, >=latex, line width=1pt] (S.north) -- (0,1.5) -- (7,1.5) node[midway,above,xshift=-1.5cm]{$\beta (1-p) S I$} -- (I.north);
\draw [->, thick, >=latex, line width=1pt] (E.north) -- (3.5,2) -- (10.5,2) node[midway,below,xshift=1.5cm]{$\eta E$} -- (R.north);
\draw [->, thick, >=latex, line width=1pt] (S.south) -| (0,-1.3) node[near end,left]{$\mu S$};
\draw [->, thick, >=latex, line width=1pt] (E.south) -| (3.5,-1.3) node[near end,left]{$\mu E$};
\draw [->, thick, >=latex, line width=1pt] (I.south) -| (7,-1.3) node[near end,left]{$\mu I$};
\draw [->, thick, >=latex, line width=1pt] (R.south) -| (10.5,-1.3) node[near end,left]{$\mu R$};
\end{tikzpicture}
\captionof{figure}{Flow diagram of the SEIR reaction-diffusion model with vital dynamics.}\label{F1}
\end{figure}
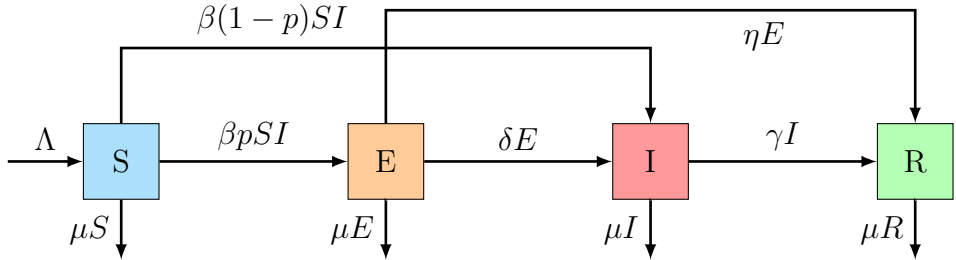

We propose a four-dimensional epidemiological model to describe the transmission dynamics of an epidemic. 
The total population is divided into four compartments: susceptible ($S$), exposed ($E$), infected ($I$), and recovered ($R$). 
Individuals are recruited into the susceptible class at a constant rate $\Lambda$ and die at a natural rate $\mu$. 
Susceptible individuals may become infected upon contact with infected individuals at a rate $\beta$, with a fraction $p$ entering the exposed class and the complementary fraction $1-p$ progressing directly to the infected class. 
Exposed individuals can either progress to the infected class at rate $\delta$ or recover without becoming infectious at rate $\eta$. Infected individuals recover at a rate $\gamma$. 
Recovered individuals acquire permanent immunity and typically only die through natural causes.
The transitions between compartments are summarized in Figure~\ref{F1}.
Spatial mobility is incorporated through diffusion terms, with a common diffusion coefficient $\lambda>0$ for all compartments. This allows the model to capture the heterogeneity in disease spread across a spatial domain. 
The model parameters are given in Table~\ref{Tab1}.
\begin{table}[H]
\centering
\setlength{\tabcolsep}{0.8cm}
\caption{Epidemic parameters for the proposed SEIR reaction-diffusion model.}\label{Tab1}
\adjustbox{max width=\textwidth}{
\begin{tabular}{cc}
\hline 
\textbf{Symbol} & \textbf{Description} \\
\hline \hline 
$\Lambda$ & Recruitment rate (immigration or birth) \\
$\mu$ & Natural mortality rate\\
$\beta$ & Transmission coefficient\\
$p$ & Fraction of infections entering the exposed class\\
$\delta$ & Progression rate from exposed to infected\\
$\eta$ & Recovery rate from the exposed class\\
$\gamma$ & Recovery rate from the infected class\\
$\lambda$ & Diffusion coefficient (identical for all classes)\\
\hline
\end{tabular}
}
\end{table}

Let $\Omega \subset \mathbb{R}^d$ be a bounded domain with smooth boundary $\partial \Omega$. The spatiotemporal dynamics of the SEIR system are described by the following system of reaction-diffusion equations
\begin{equation}\label{E2.1}
\left\{\begin{aligned}
\frac{\partial S(t, x)}{\partial t} 
 &= \lambda \Delta S(t, x) + \Lambda - \beta S(t, x) I(t, x) - \mu S(t, x), \\
\frac{\partial E(t, x)}{\partial t} &= \lambda \Delta E(t, x) + \beta p S(t, x) I(t, x) - (\delta + \eta + \mu) E(t, x), \\
\frac{\partial I(t, x)}{\partial t} &= \lambda \Delta I(t, x) + \beta (1-p) S(t, x) I(t, x) + \delta E(t, x) - (\gamma + \mu) I(t, x), \\
\frac{\partial R(t, x)}{\partial t} &= \lambda \Delta R(t, x) + \eta E(t, x) + \gamma I(t, x) - \mu R(t, x).
\end{aligned}\right. \ \text{ in } \ \mathcal{U}_T, 
\end{equation}
Here, $\mathcal{U}_T = [0, T] \times \Omega$ and $\Delta = \nabla^2$ is the Laplace operator modeling spatial diffusion.
We impose homogeneous Neumann boundary conditions to describe a closed system with no population flux across the boundary
\begin{equation}\label{E2.2}
\nabla S\cdot \vec{n} = \nabla E\cdot \vec{n} = \nabla I\cdot \vec{n} = \nabla R\cdot \vec{n} = 0, 
\ \text{ on } \Sigma_T = [0, T] \times \partial\Omega,
\end{equation}
where $\vec{n}$ is the outward unit normal to $\partial\Omega$.
The initial conditions specify the initial spatial distribution of each compartment
\begin{equation}\label{E2.3}
S(0,x)=S_0(x), \ E(0,x)=E_0(x), \ I(0,x)=I_0(x), \ R(0,x)=R_0(x), \ \text{ in } \Omega.
\end{equation}

Unlike standard SEIR models, which assume a constant population size, the inclusion of recruitment and mortality ensures that the total population 
\begin{equation*}
   N(t,x) = S(t,x) + E(t,x) + I(t,x) + R(t,x),
\end{equation*}
is bounded but not constant. We obtain for $N(t,x)$ the temporal ODE
\begin{equation}\label{eq:Ndiff}
    \frac{\partial N}{\partial t} - \lambda \Delta N = \Lambda - \mu N,
\end{equation}
i.e.\ $N$ remains bounded due to the balance between births and natural deaths, see Section~\ref{S2.2}.

\subsection{Well-Posedness Analysis}\label{S2.2}
To ensure the mathematical consistency of the proposed SEIR reaction-diffusion model, we analyze the well-posedness of the system~\eqref{E2.1}--\eqref{E2.3}. 
Our goal is to verify that the system admits a unique global solution that is bounded and nonnegative for all time. 
To this end, we rewrite the model in abstract form
\begin{equation*}
\left\{\begin{aligned}
& \partial_t \psi(t) + \mathcal{T}\psi(t) = \mathcal{F}(\psi(t)),\\ 
& \psi(0) = \psi_0,
\end{aligned}\right.
\end{equation*}
where $\psi(t) = (S,E,I,R)(t,\cdot)$, the linear operator $\mathcal{T}$ is associated with diffusion terms
\begin{equation*}
\mathcal{T}\colon \quad\begin{aligned}
&\mathcal{D}_\mathcal{T} = \Bigl\{\upsilon \in \bigl( H^2(\Omega)\bigr)^4 \mid \ \nabla \upsilon_i \cdot \vec{n} = 0, 1\le i\le 4\Bigr\} 
\subset \bigl( L^2(\Omega) \bigr)^4 
\longrightarrow \bigl( L^2(\Omega) \bigr)^4,\\
&\upsilon \longrightarrow -\lambda\Delta \upsilon,
\end{aligned} 
\end{equation*}
and $\mathcal{F}$ contains the nonlinear reaction dynamics. 
First, we establish that all solutions of the system remain uniformly bounded in finite time 
using properties of analytic $C_0$-semigroups generated by $-\Delta$ under Neumann boundary conditions \cite[pp. 230-251]{pazy2012semigroups}. 
Then, multiplying each equation by the negative part of the corresponding variable, applying integration by parts, and using Cauchy-Schwarz and Gronwall’s inequalities proves that nonnegative initial data lead to nonnegative solutions, ensuring biological feasibility. 
Furthermore, the existence and uniqueness of a global strong solution follow from standard results in functional analysis since the nonlinear function $\mathcal{F}$ is Lipschitz continuous (due to the boundedness of the solution) 
and the Laplace operator is strongly elliptic. 
Finally, the solution satisfies the regularity properties
\begin{equation*}
     \psi_i \in W^{1,2}\bigl(0, T; L^2(\Omega)\bigr) \cap L^2\bigl(0, T; H^2(\Omega)\bigr) \cap L^\infty\bigl(0, T; H^1(\Omega)\bigr),
\end{equation*}
for each component, which guarantees both smoothness and stability of the dynamics.

For further details and rigorous proofs of the well-posedness analysis, see \cite[pp. 8-11]{Zinihi2025S}. This analysis shows that the proposed system retains the essential biological characteristics of epidemic dynamics. Solutions remain nonnegative when initial data is nonnegative, and the total population, $N(t, x)$, is uniformly bounded by $\Lambda/\mu$ in the long-term regime. Additionally, including the exposed compartment alongside vital dynamics offers a biologically accurate depiction of incubation and demographic turnover. 
These properties justify using the model as a reliable testing ground for further analytical and numerical investigations, including forward and inverse PINN-based studies in spatial epidemiology.

\section{Generation of Synthetic Data via NSFD Schemes}\label{S3}
As there are no publicly available real-world epidemiological datasets with the spatial resolution required for reaction-diffusion systems, we generate synthetic data by numerically solving equations \eqref{E2.1}--\eqref{E2.3}. 
To ensure that the numerical solutions preserve essential epidemiological properties such as non-negativity, boundedness and numerical stability, we use an NSFD scheme.  

Classical finite-difference (FD) discretizations, such as the explicit/implicit Euler methods or Runge-Kutta methods, can produce negative or unbounded compartment values when applied to nonlinear epidemic PDEs, see e.g.\ \cite{maamar2024nonstandard}. 
In contrast, NSFD schemes \cite{Mickens1993, Mickens2020} are explicitly designed to preserve the qualitative features of the underlying continuous system, including the positivity of solutions, the invariance of feasible regions, and the boundedness of the total population $N$. 
These properties are essential to epidemiological models because negative population densities or unbounded growth are biologically meaningless.

\subsection{NSFD Schemes}\label{S3.1}
Using the methodology outlined in a previous work of the authors \cite{Zinihi2025NSFD}, we develop 
NSFD schemes to numerically approximate the spatiotemporal SEIR model~\eqref{E2.1}. 
Let $t_n = nk$ denote the discrete time levels with a time step size of $k > 0$. 
The time derivative is approximated by
\begin{equation*}
  \frac{\partial u(t,\cdot)}{\partial t}\bigg|_{t=t_n}
  \approx \frac{u^{n+1}_{\cdot}-u^n_{\cdot}}{\phi(k)},
\end{equation*}
where $\phi(k) > 0$ is a \textit{nonstandard denominator function} satisfying $\phi(k)\to k$ as $k\to0$.
First, we describe the one-dimensional NSFD discretization, then we extend the approach to two-dimensional spatial domains.

\subsubsection{One-Dimensional NSFD Scheme}\label{S3.1.1}
Discretize the spatial domain $\Omega=[0, L]$ into $M$ grid points with spacing $h = L/M$. For $j = 1, \dots, M$ and $n \ge0$, denote $S_j^n \approx S(t_n, x_j)$, and similarly for the other compartments $E$, $I$, and $R$. 
The term $\Delta S(t_n,x_j)$ is approximated using the following nonstandard  ('skew') discretization of the Laplacian
\begin{equation*}
    \Delta_h^2 S_{j}^{n} = \frac{S_{j+1}^n - 2S_{j}^{n+1} + S_{j-1}^n}{h^2}.
\end{equation*}
to guarantee the non-negativity of the  solution.
Consequently, the NSFD discretization of the system \eqref{E2.1} in 1D reads
\begin{equation}\label{scheme1}
\left\{\begin{aligned}
\frac{S^{n+1}_{j}-S^n_{j}}{\phi(k)} 
 &= \lambda \Delta_h^2 S_{j}^{n} + \Lambda - (\beta I^n_j + \mu) S^{n+1}_{j}, \\
\frac{E^{n+1}_{j}-E^n_{j}}{\phi(k)} 
&= \lambda \Delta_h^2 E_{j}^{n} + \beta p S^{n+1}_{j} I^n_j - (\delta + \eta + \mu) E^{n+1}_{j}, \\
\frac{I^{n+1}_{j}-I^n_{j}}{\phi(k)} 
&= \lambda \Delta_h^2 I_{j}^{n} + \beta (1-p) S^{n+1}_{j} I^n_j + \delta E^{n+1}_{j} - (\gamma + \mu) I^{n+1}_{j}, \\
\frac{R^{n+1}_{j}-R^n_{j}}{\phi(k)} 
&= \lambda \Delta_h^2 R_{j}^{n} + \eta E^{n+1}_{j} + \gamma I^{n+1}_{j} - \mu R^{n+1}_{j}.
\end{aligned}\right. 
\end{equation}
%
The NSFD scheme \eqref{scheme1} 
is reformulated in a sequentially explicit form
\begin{equation}\label{E3.1}
\left\{\begin{aligned}
S_{j}^{n+1} &= 
\frac{ S_{j}^n
+ \lambda r(k)(S_{j-1}^{n}+S_{j+1}^{n})
+ \phi(k)\Lambda }
{1 + 2\lambda r(k) + \phi(k)(\beta I_{j}^n + \mu)}, \\
E_{j}^{n+1} &= 
\frac{ E_{j}^n
+ \lambda r(k)(E_{j-1}^{n}+E_{j+1}^{n})
+ \phi(k)\beta p\, S_{j}^{n+1} I_{j}^{n} }
{1 + 2\lambda r(k) + \phi(k)(\delta+\eta+\mu)}, \\
I_j^{n+1} &= 
\frac{ I_{j}^n
+ \lambda r(k)(I_{j-1}^n + I_{j+1}^n)
+ \phi(k)\bigl[\beta(1-p) S_j^{n+1} I_j^n + \delta E_j^{n+1}\bigr] }
{1 + 2\lambda r(k) + \phi(k)(\gamma+\mu)}, \\
R_j^{n+1} &= 
\frac{ R_{j}^n
+ \lambda r(k)(R_{j-1}^{n}+R_{j+1}^{n})
+ \phi(k)[\eta E_{j}^{n+1} + \gamma I_{j}^{n+1}] }
{1 + 2\lambda r(k) + \phi(k)\mu},
\end{aligned}\right.
\end{equation}
where we used the abbreviation $r(k) =\phi(k)/h^2$ for the generalized parabolic mesh ratio.
From \eqref{scheme1} we obtain the difference equation for the discrete total population
\begin{equation}
\frac{N^{n+1}_j - N^n_j}{\phi(k)} =
\lambda \Delta_h^2 N_j^n + \Lambda - \mu N^{n+1}_j,\quad n=0,1,2,\dots, 
\end{equation}
which is the discrete analogue to \eqref{eq:Ndiff}.

\subsubsection{Two-Dimensional NSFD Scheme}\label{S3.1.2}

For a two-dimensional domain $\Omega=[0,L_x]\times[0,L_y]$ with uniform mesh spacings $h_x = h_y = h$, let $S_{j, \ell}^n \approx S(t_n, x_j, y_\ell)$ denote the pointwise approximation, with indices $j$, $\ell$ for space and $n$ for time. 
Using the same methodology as in the previous section with
\begin{equation*} 
\Delta_h^2 S_{j,l}^{n} = \frac{S_{j+1,l}^n - 2S_{j,l}^{n+1} + S_{j-1, l}^n}{h^2} + \frac{S_{j,l+1}^n - 2S_{j,l}^{n+1} + S_{j,l-1}^n}{h^2},
\end{equation*}
the 2D NSFD scheme reads
\begin{equation}\label{E3.2}
\left\{\begin{aligned}
S_{j,\ell}^{n+1} &= 
\frac{ S_{j,\ell}^n
+ \lambda r(k)\sum_{\sigma=\pm1}
(S_{j+\sigma,\ell}^n + S_{j,\ell+\sigma}^n)
+ \phi(k)\Lambda}
{1 + 4\lambda r(k)
+ \phi(k)(\beta I_{j,\ell}^n + \mu)}, \\
E_{j,\ell}^{n+1} &= 
\frac{ E_{j,\ell}^n
+ \lambda r(k)\sum_{\sigma=\pm1}
(E_{j+\sigma,\ell}^n + E_{j,\ell+\sigma}^n)
+ \phi(k)\beta p\, S_{j,\ell}^{n+1} I_{j,\ell}^{n} }
{1 + 4\lambda r(k) + \phi(k)(\delta+\eta+\mu)}, \\
I_{j,\ell}^{n+1} &= 
\frac{ I_{j,\ell}^n
+ \lambda r(k)\sum_{\sigma=\pm1}
(I_{j+\sigma,\ell}^{n} + I_{j,\ell+\sigma}^n)
+ \phi(k)\bigl[\beta(1-p) S_{j,\ell}^{n+1} I_{j,\ell}^n + \delta E_{j,\ell}^{n+1}\bigr] }
{1 + 4\lambda r(k) + \phi(k)(\gamma+\mu)}, \\
R_{j,\ell}^{n+1} &= 
\frac{ R_{j,\ell}^n
+ \lambda r(k)\sum_{\sigma=\pm1}
(R_{j+\sigma,\ell}^n + R_{j,\ell+\sigma}^n)
+ \phi(k) [ \eta E_{j,\ell}^{n+1} + \gamma I_{j,\ell}^{n+1} ] }
{1 + 4\lambda r(k) + \phi(k)\mu}.
\end{aligned}\right.
\end{equation}

\subsubsection{Qualitative Properties of the Scheme}\label{S3.1.3}
The NSFD schemes \eqref{E3.1}--\eqref{E3.2} preserve several essential qualitative properties of the continuous epidemic model. 
In particular, nonnegative initial conditions guarantee nonnegative numerical solutions for all compartments and for any choices of the 
step sizes $k,h>0$. 
Moreover, the \textit{discrete total population} $N_{j,\ell}^n = S_{j,\ell}^n + E_{j,\ell}^n + I_{j,\ell}^n + R_{j,\ell}^n$ remains uniformly bounded by $\Lambda/\mu$, consistently with the continuous formulation.
The implicit treatment of the reactive terms together with the nonstandard denominator function $\phi(k)$ also provides enhanced numerical stability, allowing the use of relatively large time steps $k$ without numerical blow-up.



\subsection{Numerical Properties and Rationale for the NSFD Scheme}\label{S3.2}
The construction of the above NSFD scheme adheres to Mickens’s foundational principles for structure-preserving discretizations \cite{Mickens1993, Mickens2020}. 
Standard FD methods often fail to maintain qualitative consistency when applied to nonlinear epidemiological PDEs, leading to negative densities or instability. 
In contrast, the present formulation introduces a \textit{nonstandard denominator function} $\phi(k)$ 
and a \textit{local nonlinear stabilizer} $\nu(U_j^n)$, that collectively ensure positivity, boundedness, and improved temporal consistency.  

This framework builds on advances in NSFD schemes for reaction-diffusion systems \cite{Zinihi2025NSFD, ChenCharpentier2013, Pasha2023, Ehrhardt2013}. 
Specifically, the local correction term $\nu(U_j^n)$ compensates for the order reduction arising from the positivity-preserving Laplacian discretization, as demonstrated in \cite{Zinihi2025NSFD}. 
The resulting scheme achieves first-order accuracy in time and second-order accuracy in space, with a local truncation error $\mathcal{O}(k) + \mathcal{O}(h^2)$.  

The positivity of the denominator function $\phi(k)$ for all $k>0$, guarantees the unconditional nonnegativity of the discrete solution.
Meanwhile, the implicit-explicit coupling of the diffusive and reactive terms enhances the numerical stability,
even for relatively large parabolic ratios $r=k/h^2$. 
These results align with stability analyses reported in \cite{ChenCharpentier2013, Pasha2023, Zinihi2025NSFD}, confirming the reliability of the NSFD approach in generating stable, high-fidelity datasets for data-driven modeling frameworks.


For both the 1D and 2D schemes, and in line with 
\cite{Zinihi2025NSFD, Ehrhardt2013}, the non-standard denominator function is chosen as the identity function $\phi(k)=k$ for the sake of simplicity, in order to simplify the proposed data generation schemes.
For a more advanced choice of $\phi(k)$ originating from the asymptotic behaviour of the total population for $t\to\infty$, we refer the reader to \cite{maamar2024nonstandard}.

\section{The Physics-Informed Neural Network Framework}\label{S4}
In this section, we describe the PINN framework that we used to approximate the solutions of the SEIR reaction-diffusion system~\eqref{E2.1}. 
PINNs integrate observational data, initial and boundary conditions, and the governing PDE system into a unified loss function.
This enables both forward simulation and inverse parameter estimation.

\subsection{The Neural Network Architecture}\label{S4.1}

We approximate the state variables
\begin{equation*}
    (S, E, I, R)(t,x) \approx  \mathcal{NN}_\theta^{\Theta_p}(t, x) := (\hat{S}_\theta, \hat{E}_\theta, \hat{I}_\theta, \hat{R}_\theta)(t,x),
\end{equation*}
where $\theta$ denotes the vector 
of trainable neural network parameters (so-called hyperparameters), including all network weights, biases, and, in the inverse problem setting, a subset of unknown epidemiological parameters to be identified from data.
Here, the set of trainable physical parameters is chosen as
\begin{equation*}
\Theta_p = (\beta,\, \delta,\, \gamma,\, \lambda)^\top\in\mathbb{R}^4.
\end{equation*}
These parameters in $\Theta_p$ were selected because they directly influence the dynamics of epidemic transmission, the temporal progression of the disease, and the mechanisms of spatial spreading. 
Furthermore, they are typically challenging to measure accurately and can vary significantly across outbreaks, geographical regions, and intervention strategies. 
Thus, estimating these parameters from observational data is an important inverse problem in epidemiological modeling.

In contrast, the demographic parameters $\Lambda$ and $\mu$ are assumed to be known constants and are not included in the training procedure. 
These quantities correspond to long-term demographic processes, such as birth rates and natural mortality, which evolve on slower timescales than epidemic dynamics. 
They can often be estimated independently using 
demographic statistics.

Furthermore, the parameter $p$ and the transition coefficient $\eta$ 
are treated as fixed quantities to reduce the dimensionality of the inverse problem and improve parameter identifiability during training. 
Estimating too many parameters simultaneously can cause optimization instability 
and non-uniqueness issues, especially when observational data is sparse or noisy. 
Thus, the chosen subset of parameters in $\Theta_p$ strikes a balance between epidemiological relevance, practical identifiability, and computational stability within the PINN framework.
%
Each compartment is approximated by a fully connected feedforward neural network.
The unknown parameters in $\Theta_p$ are optimized jointly with the neural network weights using gradient-based optimization and automatic differentiation.

The inputs to the neural network are the temporal coordinate $t \in [0, T]$ and the spatial coordinates $x\in\Omega$ for $d=1$, or $(x,y)\in\Omega$ for $d=2$, while the outputs correspond to the four epidemiological compartments $(\hat{S},\hat{E},\hat{I},\hat{R})$.
To improve the representation of high-frequency spatial patterns and mitigate the spectral bias commonly encountered in standard multilayer perceptrons, \textit{Fourier feature mappings} \cite{Wang2021S, Tancik2020} are employed for the spatial inputs. 
Specifically, the spatial coordinates are transformed according to
\begin{equation*}
    \gamma(x) = \bigl[\sin(2\pi \mathbf{B}x),\cos(2\pi \mathbf{B}x) \bigr],
\end{equation*}
where $\mathbf{B}$) is a matrix of randomly sampled frequencies. 
This encoding significantly enhances the expressive power of the PINN, enabling the accurate reconstruction of sharp epidemic fronts and localized infection waves in both one- and two-dimensional domains. 
The network consists of $L$ hidden layers with nonlinear activations (\texttt{tanh} or \texttt{swish}), and residual skip connections are incorporated to improve trainability and stability for long-time dynamics.

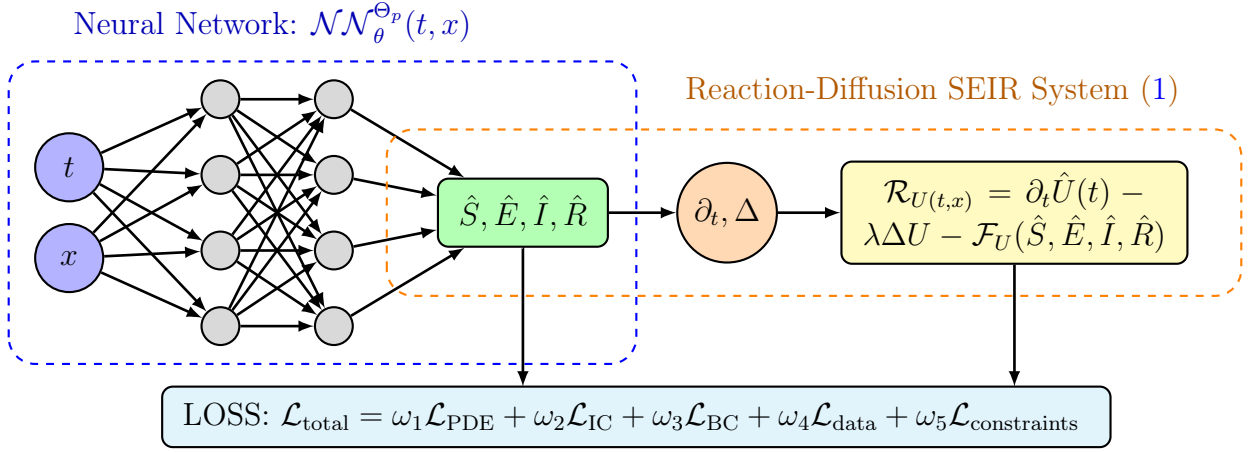
\begin{figure}[H]
\centering
\begin{tikzpicture}
\draw[rounded corners=10pt, dashed, thick, blue] (-4.8, -1) rectangle (3.5, 3) {};
\node[circle, draw, fill=blue!30, minimum size=0.9cm, thick, xshift=-1cm, yshift=1cm] (inputt) at (-3, 0.6) {$t$};
\node[circle, draw, fill=blue!30, minimum size=0.9cm, thick, xshift=-1cm, yshift=1cm] (inputx) at (-3, -0.6) {$x$};
\node[circle, draw, fill=gray!30, minimum size=0.5cm, thick, yshift=1cm] (h1) at (-2, 1.5) {};
\node[circle, draw, fill=gray!30, minimum size=0.5cm, thick, yshift=1cm] (h2) at (-2, 0.5) {};
\node[circle, draw, fill=gray!30, minimum size=0.5cm, thick, yshift=1cm] (h3) at (-2, -0.5) {};
\node[circle, draw, fill=gray!30, minimum size=0.5cm, thick, yshift=1cm] (h4) at (-2, -1.5) {};
\node[circle, draw, fill=gray!30, minimum size=0.5cm, thick, yshift=1cm] (h5) at (-0.5, 1.5) {};
\node[circle, draw, fill=gray!30, minimum size=0.5cm, thick, yshift=1cm] (h6) at (-0.5, 0.5) {};
\node[circle, draw, fill=gray!30, minimum size=0.5cm, thick, yshift=1cm] (h7) at (-0.5, -0.5) {};
\node[circle, draw, fill=gray!30, minimum size=0.5cm, thick, yshift=1cm] (h8) at (-0.5, -1.5) {};
\node[rectangle, rounded corners, draw, fill=green!30, minimum size=0.9cm, thick, xshift=0.5cm, yshift=1cm] (output) at (1.5, 0) { \ $\hat{S}, \hat{E}, \hat{I}, \hat{R}$ \ };
\foreach \i in {h1,h2,h3,h4}
  \draw[->, thick, >=latex, line width=1pt] (inputt) -- (\i);
\foreach \i in {h1,h2,h3,h4}
  \draw[->, thick, >=latex, line width=1pt] (inputx) -- (\i);
\foreach \i in {h5,h6,h7,h8}
  \foreach \j in {h1,h2,h3,h4}
    \draw[->, thick, >=latex, line width=1pt] (\j) -- (\i);
\foreach \i in {h5,h6,h7,h8}
  \draw[->, thick, >=latex, line width=1pt] (\i) -- (output);
\node[anchor=north,above] at (-0.9, 3.1) {\color{blue!70!black}Neural Network: $\mathcal{NN}_\theta^{\Theta_p}(t, x) \qquad$};
\draw[rounded corners=10pt, dashed, orange, thick, yshift=0.9cm] (0.2, -1) rectangle (11.5, 1.2) {};
\node[circle, draw, fill=orange!30, thick, yshift=1cm] (dudt) at (4.7, 0) {$\partial_t, \Delta$};
\node[rectangle, rounded corners, draw, fill=yellow!30, thick, align=center, yshift=1cm, text width=4.3cm] (residual) at (8.5, 0) {$\mathcal{R}_{U(t, x)} = \partial_t \hat{U}(t) - \lambda\Delta U - \mathcal{F}_U(\hat{S}, \hat{E}, \hat{I}, \hat{R})$};
\draw[->, thick, >=latex, line width=1pt] (output) -- (dudt);
\draw[->, thick, >=latex, line width=1pt] (dudt) -- (residual);
\node[anchor=north, above, yshift=1cm] at (6.25, 1.25) { \hspace{2.2cm} \color{orange!70!black}Reaction-Diffusion SEIR System~\eqref{E2.1}};
\node[rectangle, rounded corners, draw, fill=cyan!10, thick, align=center, minimum size=0.8cm, xshift=-0.5cm, yshift=0.8cm] (loss) at (4, -2.5) { \ LOSS: $\mathcal{L}_{\text{total}} = \omega_1 \mathcal{L}_{\text{PDE}} + \omega_2 \mathcal{L}_{\text{IC}} + \omega_3 \mathcal{L}_{\text{BC}} + \omega_4 \mathcal{L}_{\text{data}} + \omega_5 \mathcal{L}_{\text{constraints}}$ \ };
\draw [->, thick, >=latex, line width=1pt] (output.south) -| (2,-1.3);
\draw [->, thick, >=latex, line width=1pt] (residual.south) -| (8.5,0) -- (8.5,-1.3);
\end{tikzpicture}
\captionof{figure}{Schematic diagram of the proposed PINN architecture for the reaction-diffusion SEIR model~\eqref{E2.1}--\eqref{E2.3}.}\label{F2}
\end{figure}

    
Figure~\ref{F2} illustrates the architecture of the proposed SEIR–PINN. 
The neural network takes the time $t$ and the spatial coordinates $x\in\Omega$ as inputs and outputs the approximate compartmental fields $(\hat{S}, \hat{E}, \hat{I}, \hat{R})$. 
Automatic differentiation computes the temporal and spatial derivatives, enabling the evaluation of the PDE residuals corresponding to the reaction-diffusion SEIR system~\eqref{E2.1}. 
Together with the data misfit, boundary and initial condition losses, and epidemiological constraints, these residuals contribute to the total composite loss~\eqref{E4.1}, which guides the training process.
These residuals, along with the data misfit, boundary and initial condition losses, and epidemiological constraints, contribute to the total composite loss~\eqref{E4.1}, which guides the training process.

\subsection{The Loss Function Design}\label{S4.2}

The total loss function is defined as a weighted sum of multiple components
\begin{equation}\label{E4.1}
    \mathcal{L}_{\text{total}} = \omega_1 \mathcal{L}_{\text{PDE}} + \omega_2 \mathcal{L}_{\text{IC}} + \omega_3 \mathcal{L}_{\text{BC}} + \omega_4 \mathcal{L}_{\text{data}} + \omega_5 \mathcal{L}_{\text{constraints}},
\end{equation}
where each component is explained below.
\paragraph{(i) PDE Residual Loss}  
For each compartment, the PDE residual is computed via automatic differentiation (AD) of the neural outputs
\begin{equation*}
      \mathcal{R}_S = \partial_t \hat{S} - \lambda \Delta \hat{S} - \Lambda + \beta \hat{S}\hat{I} + \mu \hat{S},
\end{equation*}
with analogous definitions $\mathcal{R}_E$, $\mathcal{R}_I$, and $\mathcal{R}_R$. 
The PDE loss is
\begin{equation*}
    \mathcal{L}_{\text{PDE}} = \frac{1}{N_r}\sum_{i=1}^{N_r} \bigl( \mathcal{R}_S^2 + \mathcal{R}_E^2 + \mathcal{R}_I^2 + \mathcal{R}_R^2 \bigr)(t_i, x_i),
\end{equation*}
evaluated at $N_r$ randomly sampled collocation points $(t_i, x_i)\in \mathcal{U}_T$.

\paragraph{(ii) Initial Condition Loss}  
Enforces agreement with prescribed initial data
\begin{equation*}
\mathcal{L}_{\text{IC}} = \frac{1}{N_{\text{IC}}}\sum_{i=1}^{N_{\text{IC}}} \Big( 
\bigl|\hat{S}_0(x_i) - S_0(x_i)\bigr|^2 +
\bigl|\hat{E}_0(x_i) - E_0(x_i)\bigr|^2 +
\bigl|\hat{I}_0(x_i) - I_0(x_i)\bigr|^2 +
\bigl|\hat{R}_0(x_i) - R_0(x_i)\bigr|^2
\Big),
\end{equation*}
where $N_{IC}$ denotes the number of sampled points 
at 
$t = 0$.

\paragraph{(iii) Boundary Condition Loss}  
Homogeneous Neumann conditions are enforced weakly by penalizing the normal derivative
\begin{equation*}
\mathcal{L}_{\text{BC}} = \frac{1}{N_{\partial\Omega}}\sum_{i=1}^{N_{\partial\Omega}}
\Bigl( \bigl|\nabla \hat{S}\cdot \vec{n}\bigr|^2
+ \bigl|\nabla \hat{E}\cdot \vec{n}\bigr|^2 
+ \bigl|\nabla \hat{I}\cdot \vec{n}\bigr|^2 
+ \bigl|\nabla \hat{R}\cdot \vec{n}\bigr|^2 \Bigr)(t_i, x_i),
\end{equation*}
for $N_{\partial\Omega}$ sampled boundary points $(t_i, x_i)\in \Sigma_T$.

\paragraph{(iv) Data Loss}  
When synthetic observations are available at discrete time-space locations $\{(t_j,x_j)\}$, 
we minimize the misfit between the network outputs and the data values $\mathcal{D}_j$
\begin{equation*}
\mathcal{L}_{\text{data}} = \frac{1}{N_d} \sum_{j=1}^{N_d} \bigl| \hat{U}(t_j,x_j) - \mathcal{D}_j \bigr|^2,
\end{equation*}
where $N_d$ represents the number of observational data points and $U$ denotes the observed compartments (e.g., $I$ or $(I,R)$).

\paragraph{(v) Constraint Loss}  
To ensure epidemiological consistency, the PINN framework incorporates additional constraint-based penalty terms within the loss function. 
Non-negativity of the compartment populations is enforced through soft penalties of the form
$\max(0,-\hat{U})^2$ penalties for $U\in\{S,E,I,R\}$.
Furthermore, deviations from the expected population balance induced by the vital dynamics are penalized via
\begin{equation*}
    \mathcal{L}_{\text{pop}} = \frac{1}{N_r}\sum_{i=1}^{N_r} \Big( \hat{S}+\hat{E}+\hat{I}+\hat{R} - \frac{\Lambda}{\mu} \Big)_+^2,
\end{equation*}
where $(\cdot)_+$ denotes the positive part, ensuring that the total population does not exceed its carrying capacity. 
In addition, admissibility of learned epidemiological parameters is guaranteed through bounded transformations such as
 $\beta = \beta_{\max}\,\sigma(\tilde{\beta})$,
 where $\sigma$ denotes the logistic sigmoid function, ensuring $0<\beta<\beta_{\max}$.
 Altogether, the constraint contribution to the optimization problem is given by
 \begin{equation*}
\mathcal{L}_{\text{constraints}} = \mathcal{L}_{\text{nonneg}} + \mathcal{L}_{\text{pop}} + \mathcal{L}_{\text{param}}.
\end{equation*}




In the present work, we formulate the components of the composite loss function using the standard $L^2$-based mean squared error (MSE). The $L^2$-based MSE is widely used in the PINN literature due to its smooth optimization properties and compatibility with gradient-based training methods. 
Recent studies have shown that alternative losses, such as $L^1$-based formulations, may improve robustness in certain situations \cite{wang2022}. 
However, we use the $L^2$-loss in this work because the considered SEIR reaction-diffusion system produces sufficiently smooth solutions and stable residual distributions. This makes the $L^2$ formulation well-suited for our framework.

\subsection{The Training Strategy}\label{S4.3}
The training procedure integrates the composite loss function defined in Section~\ref{S4.2} with a staged optimization strategy. 
This strategy is designed to enhance convergence stability and parameter identifiability. The complete workflow is summarized below in Algorithm~\ref{Alg1}.

\begin{algorithm}[htb]
\caption{SEIR-PINN Training Algorithm}\label{Alg1}
\begin{algorithmic}[1]
\STATE \emph{Initialization:} 
Neural network parameters $\theta$ are initialized using Xavier initialization~\cite{Glorot2010}. 
When the inverse problem is considered, the epidemiological and diffusion parameters $\Theta_p = (\beta, p, \delta, \eta, \gamma, \lambda)$ are initialized within biologically admissible ranges inferred from prior studies, while demographic quantities $(\Lambda, \mu)$ are fixed.
\STATE \emph{Staged Optimization:}
Training is performed in three successive stages to promote robustness and accuracy:
\begin{itemize}\setlength{\itemsep}{0cm}
\item[$i:$] \emph{Stage I (Data fitting):} 
The network is pretrained by minimizing $\mathcal{L}_{\text{data}} + \mathcal{L}_{\text{IC}}$ using the Adam optimizer~\cite{Liu1989} to match observations and initial conditions.

\item[$ii:$] \emph{Stage II (Physics enforcement):} 
The PDE residual and boundary terms are introduced, and the total composite loss $\mathcal{L}_{\text{total}}$ in~\eqref{E4.1} is minimized using Adam with a reduced learning rate to ensure the stable incorporation of physical constraints.

\item[$iii:$] \emph{Stage III (Refinement):} 
The optimization switches to the quasi-Newton L-BFGS method~\cite{Zinihi2025PINN, Raissi2019, Liu1989} for fine-tuning until convergence, achieving high-precision satisfaction of the governing equations.
\end{itemize}
\STATE \emph{Adaptive Sampling:}  
To improve sample efficiency and numerical accuracy, collocation points are resampled dynamically at each epoch, with a higher probability in regions with large PDE residuals. 
Training is accelerated for two-dimensional domains via GPU parallelization and mixed-precision arithmetic.
\STATE \emph{Monitoring and Constraints:}  
Non-negativity, boundedness, and population balance constraints (Section~\ref{S4.2}) are softly enforced during training. 
Throughout the process, the total loss 
and its individual components are monitored, and early stopping is applied when the relative improvement falls below a predefined tolerance threshold.
\end{algorithmic}
\end{algorithm}

The proposed SEIR-PINN framework enables three complementary learning modes:\\
\emph{Forward simulation:} Approximating the full spatiotemporal evolution $(S,E,I,R)(t,x)$ given initial and boundary conditions.\\
\emph{Inverse identification:} Estimating unknown epidemiological and diffusion parameters from sparse or noisy data.\\
\emph{Constraint-aware prediction:} Ensuring epidemiological plausibility by embedding population conservation and non-negativity directly into the optimization.

Algorithm~\ref{Alg1} thus provides a stable and interpretable training process capable of reconstructing epidemic dynamics and inferring key transmission parameters with high fidelity.

\section{Numerical Results}\label{S5}

In this section, we present and discuss the numerical results obtained from Sections~\ref{S3} and~\ref{S4} applied to the one- and two-dimensional SEIR reaction-diffusion system~\eqref{E2.1}--\eqref{E2.3}. 
The PINN is trained to simultaneously solve the forward problem and recover the unknown epidemiological parameters $\{\beta,\delta,\gamma,\lambda\}$ via the inverse problem formulation described in Section~\ref{S4}.

Synthetic ground-truth data are generated using the NSFD schemes, namely \eqref{E3.1} in the one-dimensional setting and \eqref{E3.2} in the two-dimensional setting, with both formulations enforcing zero-flux Neumann boundary conditions.

We initialize the system with normalized initial conditions, where each compartment is expressed as a proportion of the total population. 
All experiments are conducted using the parameter values listed in Table~\ref{Tab2}, and the true values of the inverse parameters are used as reference targets for evaluating the recovery accuracy. Parameters marked with $(\ast)$ are estimated through the PINN inverse problem, while the remaining parameters are kept fixed.

\begin{table}[H]
\centering
\setlength{\tabcolsep}{0.7cm}
\caption{Parameter values used in all numerical experiments.}\label{Tab2}
\adjustbox{max width=\textwidth}{
\begin{tabular}{lccccc}
\hline
\textbf{Parameter} & \textbf{Symbol} & \textbf{Value} & \textbf{Reference} &
\textbf{Role} & \textbf{Status} \\
\hline\hline
Recruitment rate & $\Lambda$ & 1.0000 & \cite{Zinihi2025NSFD, Liu2022} & Fixed & ---\\ 
Death rate & $\mu$ & 0.0100 & \cite{CDCCOVID, CDCFlu, Zinihi2026Koopman} & Fixed & --- \\
Transmission rate        & $\beta$ & 0.4000 & \cite{UKGov2025, WHO2024Covid, JHU2020, CDCCOVID, CDCFlu, Zinihi2026Koopman, ECDCFlu} & Inverse & $(\ast)$ \\
Direct exposure fraction & $p$ & 0.3000 & Assumed & Fixed & --- \\
Progression rate ($E\to I$) & $\delta$ & 0.3000 & \cite{UKGov2025, WHO2024Covid} & Inverse & $(\ast)$ \\
Recovery rate ($E\to R$) & $\eta$    & 0.1000 & \cite{WHO2024Covid, JHU2020, Zinihi2026Koopman, UKGov2025, CDCCOVID, CDCFlu} & Fixed & --- \\
Recovery rate ($I\to R$)        & $\gamma$  & 0.2000 & \cite{CDCFlu, Zinihi2026Koopman} & Inverse & $(\ast)$ \\
Diffusion coefficient    & $\lambda$ & 0.0500 & \cite{Zinihi2025NSFD, Liu2022} & Inverse & $(\ast)$ \\
\hline
\end{tabular}
}
\end{table}

\subsection{One-Dimensional Results}\label{S5.1}

We begin with the one-dimensional configuration in which the spatial domain is $\Omega = [0, L]$ with $L = 1$ (\textit{dimensionless}) and the time horizon is $T = 5$ (\textit{years}). 
The initial conditions place a Gaussian seed of exposed and infected individuals centred at $x = 0.5$, while the susceptible population is nearly uniform at $S_0 \approx 90\,\%$; the recovered compartment is initialized to zero.


The model is trained for 8000 epochs using the Adam optimizer, with an initial learning rate of $10^{-3}$ and cosine annealing scheduling down to $10^{-5}$. 
To stabilize training, gradient clipping with a maximum norm of 1.0 is employed.

\subsubsection{Space--Time Dynamics Comparison}\label{S5.1.1}

\begin{figure}[H]
\centering
\includegraphics[width=1\textwidth]{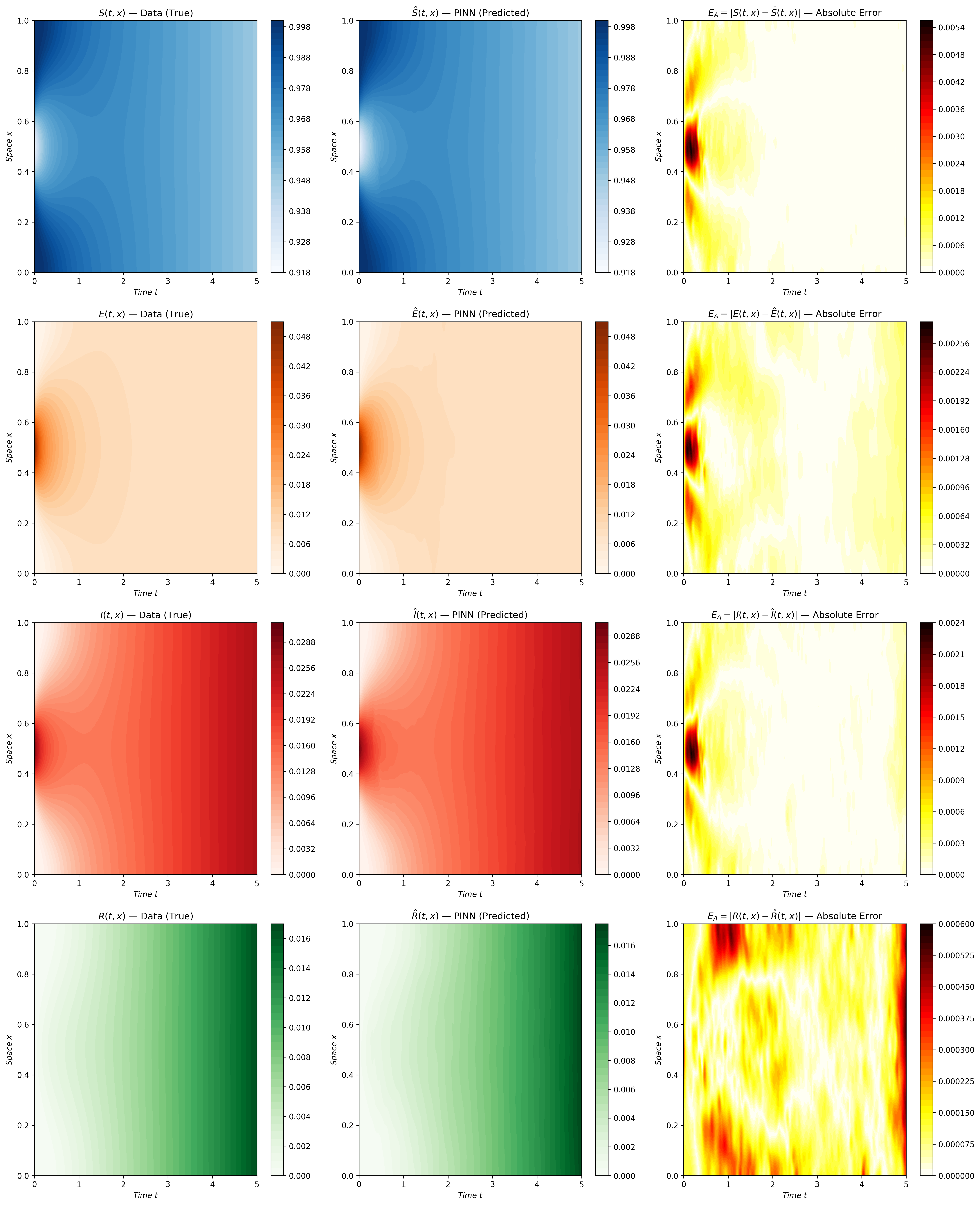}
\caption{Space–time dynamics: Reference data vs PINN prediction; 1D case.}\label{F3}
\end{figure}

Figure~\ref{F3} displays the space--time contour maps for all compartments across the domain $\mathcal{U}$, comparing the reference data (true solution, left column), the PINN prediction (centre column), and the pointwise absolute error $E_A = |\hat{U} - U|$ (right column). 
The susceptible population decreases monotonically from its initial high value as the disease spreads outward from the center due to diffusion, while the exposed and infected compartments exhibit a wave-like propagation pattern consistent with reaction–diffusion dynamics.
The absolute errors are mainly concentrated near zero, which supports that the PINN accurately reproduces the qualitative and quantitative behavior of all four compartments.

\subsubsection{Fixed-Point Temporal Profiles}\label{S5.1.2}

\begin{figure}[H]
\centering
\includegraphics[width=1\textwidth]{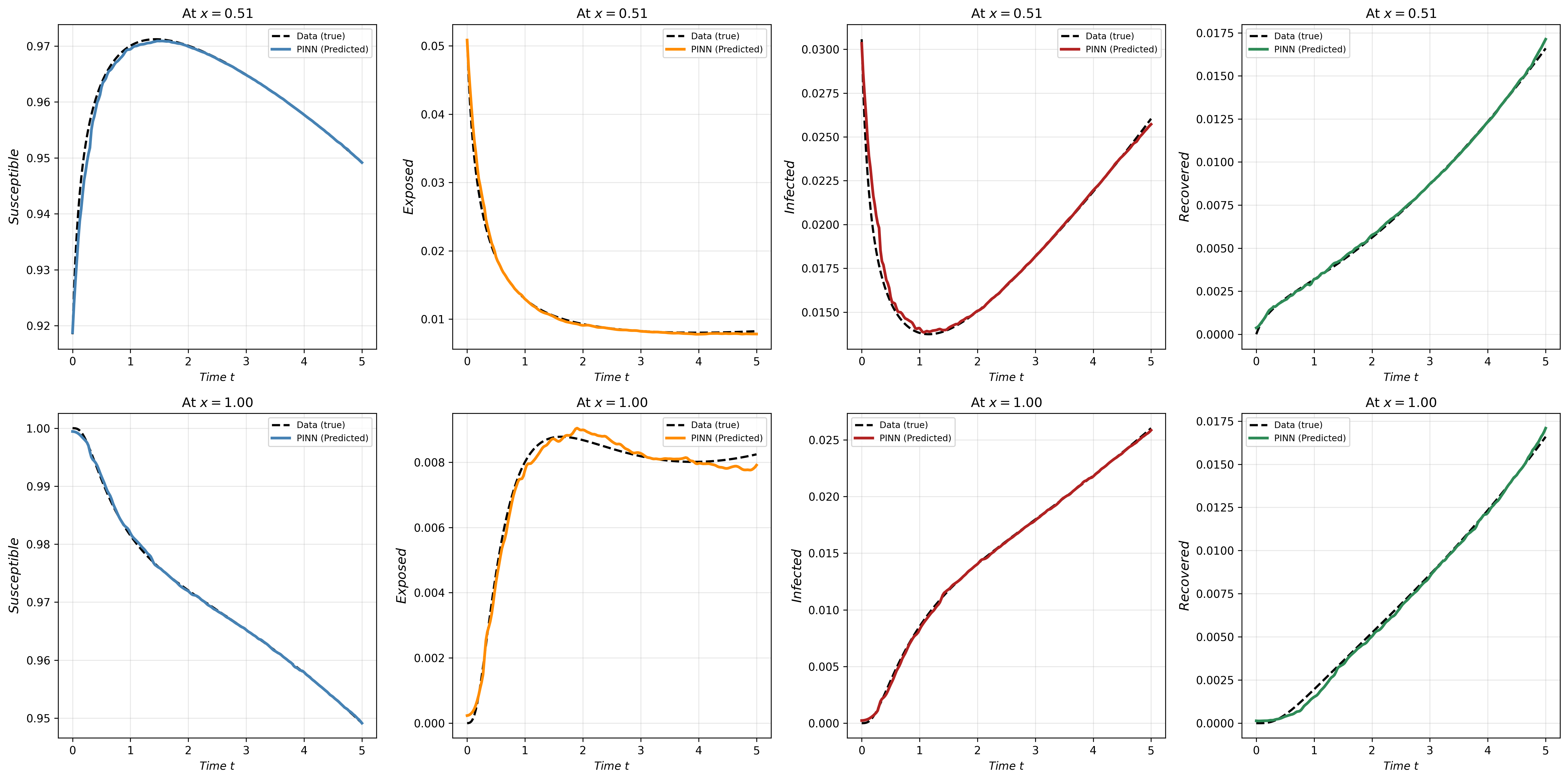}
\caption{Temporal dynamics at fixed spatial locations: Reference data vs PINN prediction; 1D case.}\label{F4}
\end{figure}

Figure~\ref{F4} compares the PINN prediction against the reference data at two fixed spatial locations: the midpoint $x = 0.50$ and the boundary point $x = 1.00$. 
In both cases, the PINN accurately reproduces the timing and magnitude of the epidemic peaks, with minor discrepancies observed for the exposed compartment at the boundary location. 
Overall, the errors remain small and are concentrated near zero.

\subsubsection{Training Loss Convergence}\label{S5.1.3}

\begin{figure}[H]
\centering
\includegraphics[width=1\textwidth]{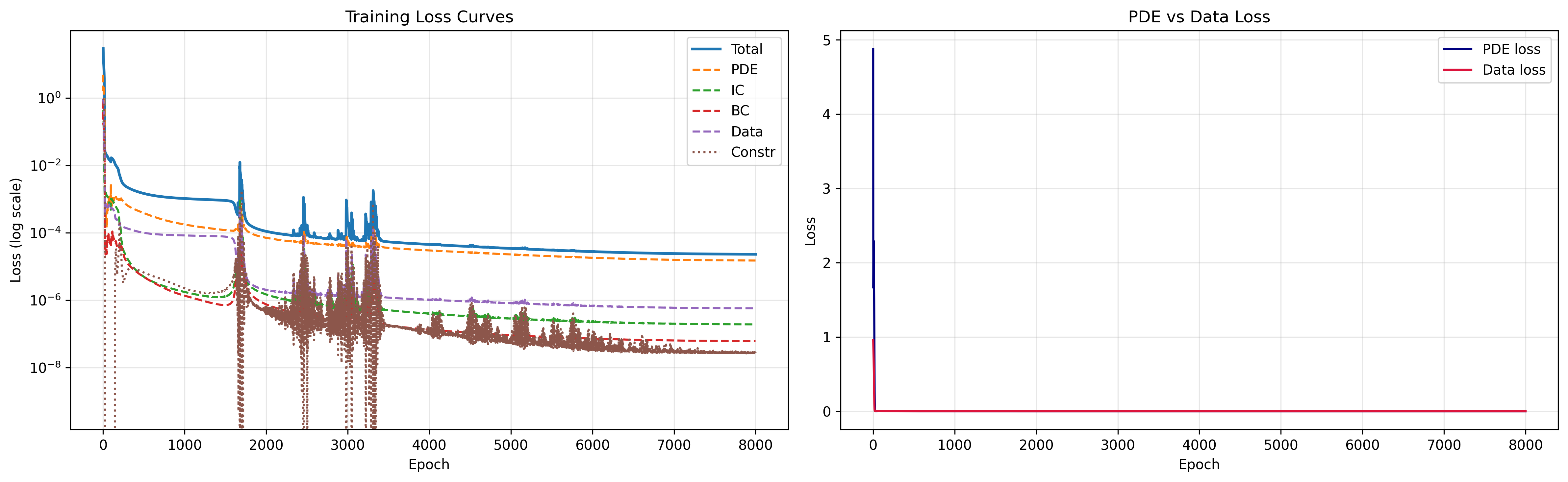}
\caption{PINN training loss evolution over 8000 epochs; 1D case.}\label{F5}
\end{figure}

Figure~\ref{F5} illustrates the evolution of each loss component over 8000 training epochs. 
The total loss decreases by approximately two to three orders of magnitude. 
The data loss and initial condition loss dominate in the early epochs, guiding the network toward the correct solution manifold.

\subsubsection{Parameter Recovery}\label{S5.1.4}

\begin{figure}[H]
\centering
\includegraphics[width=1\textwidth]{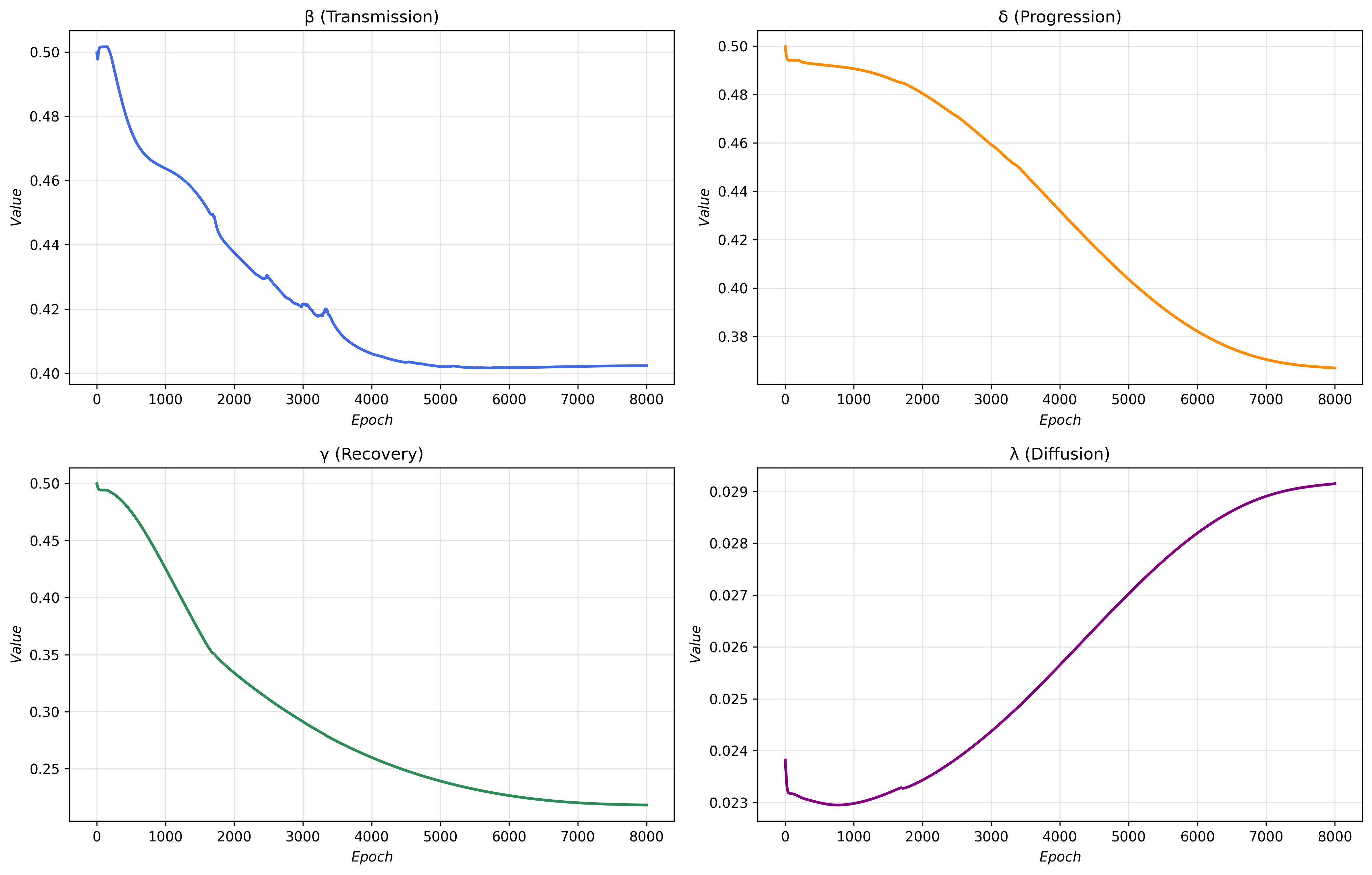}
\caption{Convergence of estimated parameters during PINN training; 1D case.}\label{F6}
\end{figure}

Figure~\ref{F6} complements this analysis by illustrating the convergence history of the four estimated parameters $\beta$, $\delta$, $\gamma$, and $\lambda$. 

\subsubsection{Quantitative Error Analysis}\label{S5.1.5}

Tables~\ref{Tab4} and~\ref{Tab5} summarize the pointwise error metrics for each compartment and the parameter recovery results, respectively.

{
\begin{table}[H]
\centering
\setlength{\tabcolsep}{0.7cm}
\caption{Quantitative error metrics between the PINN predictions and the reference data; 1D case.}\label{Tab4}
\adjustbox{max width=\textwidth}{
\begin{tabular}{ccccc}
\hline
\textbf{Compartment} & \textbf{Rel.\ $L^2$ Error} &
\textbf{MAE} & \textbf{RMSE} & \textbf{Max Error} \\
\hline\hline
$S$ & 6.0606\,e-04 & 2.5171\,e-04 & 5.8601\,e-04 & 5.5081\,e-03 \\
$E$ & 3.5408\,e-02 & 2.2358\,e-04 & 3.5202\,e-04 & 2.7316\,e-03 \\
$I$ & 1.4550\,e-02 & 1.2322\,e-04 & 2.5424\,e-04 & 2.3538\,e-03 \\
$R$ & 1.6345\,e-02 & 1.0933\,e-04 & 1.4441\,e-04 & 5.9132\,e-04 \\
\hline
\end{tabular}
}
\end{table}
\vspace{-0.8em}
\begin{table}[H]
\centering
\setlength{\tabcolsep}{0.7cm}
\caption{Comparison between the true and PINN-estimated values of the inferred parameters; 1D case.}\label{Tab5}
\adjustbox{max width=\textwidth}{
\begin{tabular}{ccccc}
\hline
\textbf{Parameter} & \textbf{True Value} & \textbf{Estimated} &
\textbf{Abs.\ Error} & \textbf{Rel.\ Error} \\
\hline\hline
$\beta$   & 0.40 & 0.4024 & 5.9132\,e-04 & 0.59\% \\
$\delta$  & 0.30 & 0.3669 & 6.6900\,e-02 & 22.30\% \\
$\gamma$  & 0.20 & 0.2183 & 1.8300\,e-02 & 9.15\% \\
$\lambda$ & 0.05 & 0.0292 & 2.0800\,e-02 & 41.60\% \\
\hline
\end{tabular}
}
\end{table}
}

\subsection{Two-Dimensional Results}\label{S5.2}

We now extend the numerical investigation to the two-dimensional setting with spatial domain $\Omega = [0, L_x] \times [0, L_y]$, 
$L_x = L_y = L$, and the same time horizon $T$. 
The added spatial dimension substantially increases both the complexity of the forward problem and the difficulty of the inverse identification task, as the Laplacian now acts on two spatial directions simultaneously and the boundary condition must be enforced on all four walls of the domain.


Training proceeds for 10,000 epochs with the same Adam optimizer and cosine
annealing schedule. 
The four Neumann boundary conditions are enforced independently for each wall: the left and right walls enforce $\partial_x \hat{U} = 0$ while the top and bottom walls enforce $\partial_y \hat{U} = 0$, with each wall contributing an independent term to the total boundary loss.

\subsubsection{Spatial Snapshots at Five Time Instances}\label{S5.2.1}

\begin{figure}[H]
 \centering
\includegraphics[width=1\textwidth]{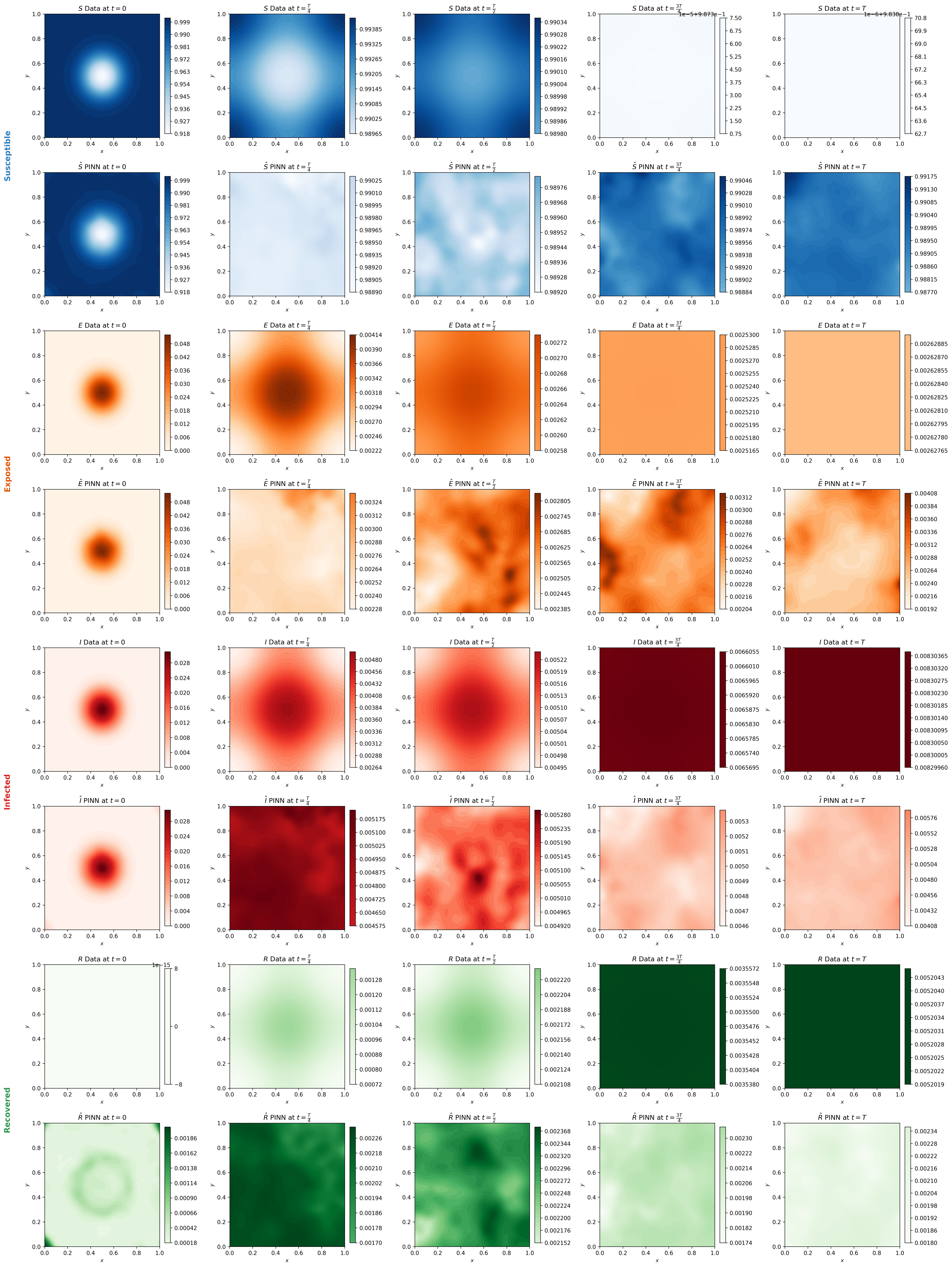}
 \caption{Spatial snapshots: Reference data vs PINN prediction; 2D case.}\label{F7}
 \end{figure}

Figure~\ref{F7} presents the spatial distribution of each compartment at five representative time instances $t = 0,\ T/4,\ T/2,\ 3T/4,\ T$. 
For each compartment, the reference data and the PINN prediction are shown as filled contour plots. 
At $t = 0$, the infectious seed is confined to a small neighborhood around the center of the domain, and both the reference data and the PINN prediction show close agreement, as indicated by the colormap scales. 
As time progresses, the disease diffuses radially outward and simultaneously grows due to the nonlinear reaction kinetics, producing a spreading wave visible in the $I$ and $E$ maps. 
The susceptible population $S$ is depleted most rapidly at the centre, while the recovered compartment $R$ accumulates from the centre outward. 
The PINN accurately reproduces this spatio-temporal evolution across all five snapshots, capturing both the diffusive spreading and the reaction-driven dynamics.

\subsubsection{Three-Dimensional Surface Comparison}\label{S5.2.2}

\begin{figure}[H]
 \centering
\includegraphics[width=1\textwidth, height=1.4\textwidth]{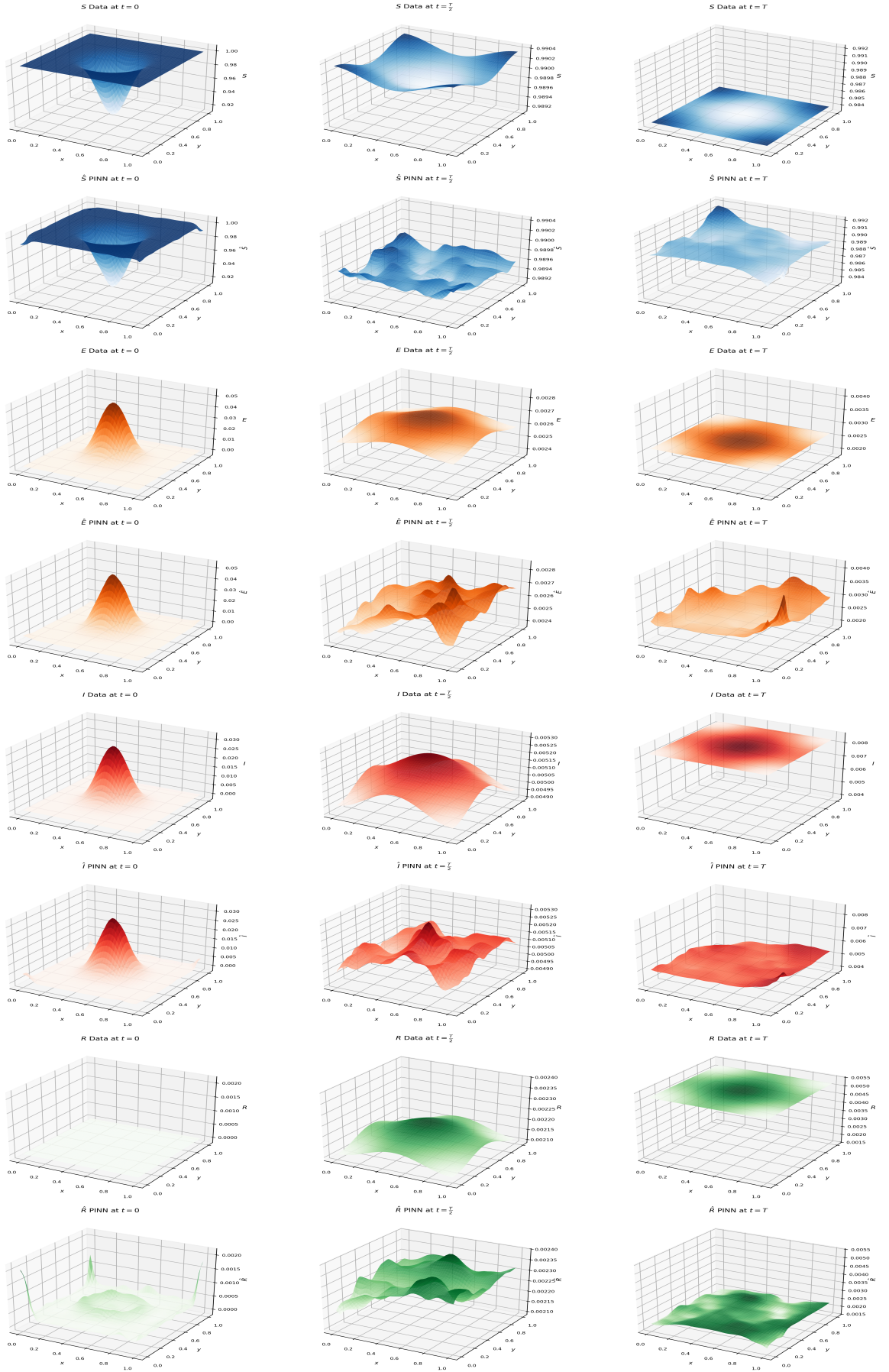}
 \caption{Three-dimensional surface: Reference data vs PINN prediction; 2D case.}\label{F8}
 \end{figure}

Figure~\ref{F8} presents 3D surface plots of the reference data and PINN prediction for all four compartments at the representative time instances $t = 0, T/2, T$. 
The surface visualizations clearly highlight the spatial heterogeneity of the compartment distributions.
The PINN surfaces closely follow the reference data in terms of height, shape, and boundary behavior, with only minor discrepancies, as further illustrated in Figure~\ref{F9} and Table~\ref{Tab7}. No spurious oscillations are observed near the domain boundaries, confirming the effectiveness of the Neumann boundary loss in the two-dimensional setting.

\subsubsection{Pointwise Error Maps}\label{S5.2.3}

\begin{figure}[H]
 \centering
\includegraphics[width=1\textwidth]{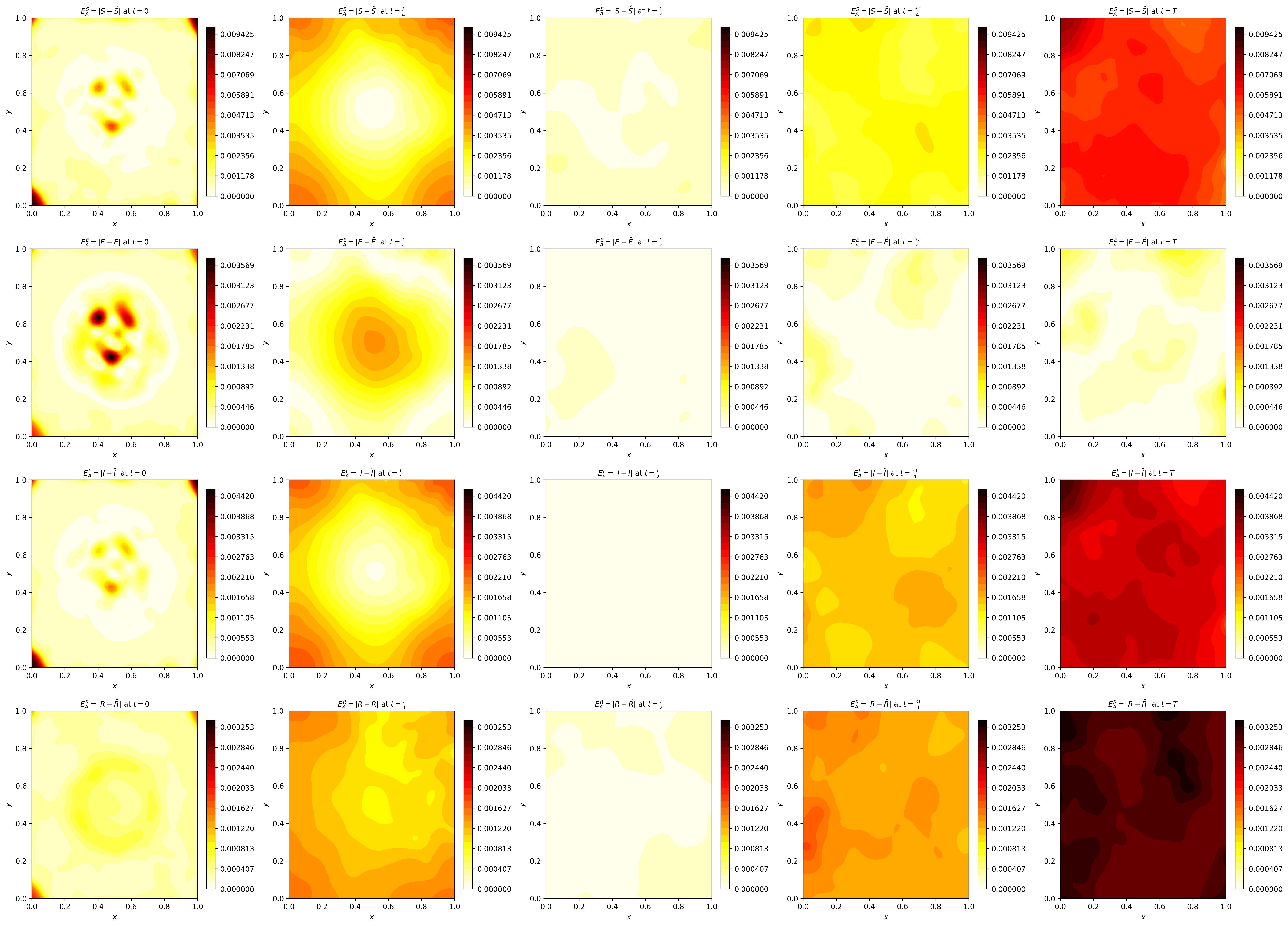}
 \caption{Pointwise error maps at five time snapshots; 2D case.}\label{F9}
 \end{figure}

Figure~\ref{F9} displays the absolute error maps $E_A = |\hat{U} - U|$ for all four compartments at five time snapshots  $t = 0,\ T/4,\ T/2,\ 3T/4,\ T$.
Importantly, the errors are distributed across the domain without any systematic spatial bias, confirming that the collocation strategy provides adequate coverage of the spatio-temporal domain.
Table~\ref{Tab7} reports the corresponding quantitative error metrics between the PINN predictions and the reference data.

\subsubsection{Training Loss Convergence}\label{S5.2.4}

\begin{figure}[H]
 \centering
\includegraphics[width=1\textwidth]{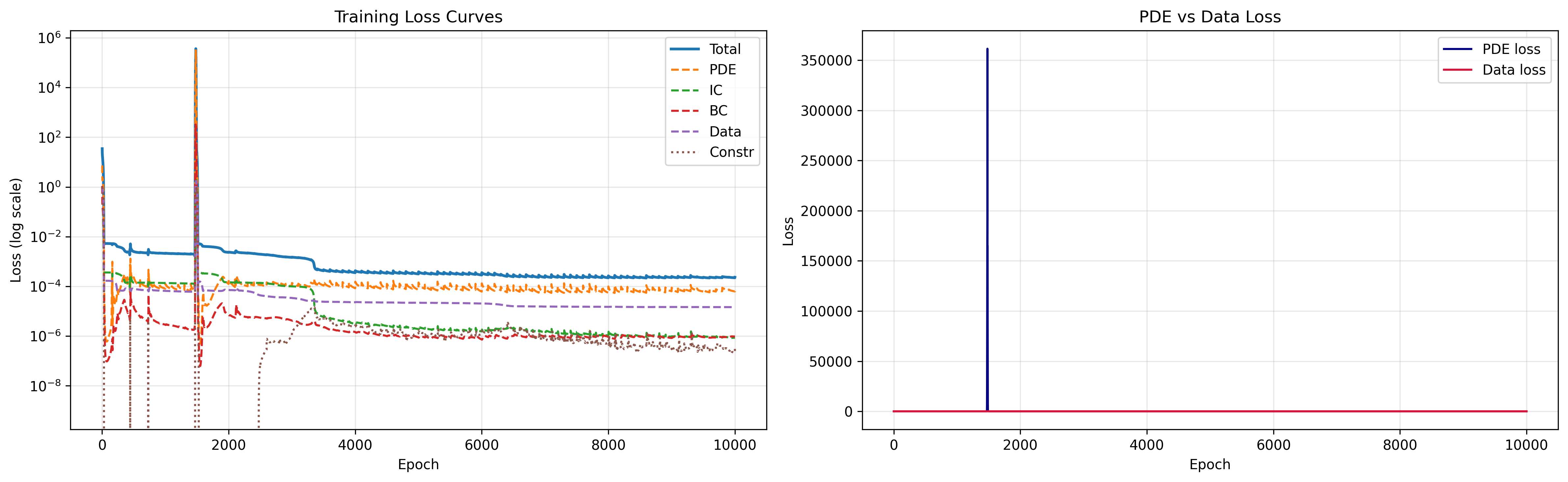}
 \caption{PINN training loss evolution over 10000 epochs; 2D case.}\label{F10}
 \end{figure}

Figure~\ref{F10} shows the loss component histories over the 10\,000 training epochs.
The overall convergence behaviour is qualitatively similar to the 1D case, Figure~\ref{F5}, but with a slower rate of decrease reflecting the greater approximation difficulty. 

\subsubsection{Parameter Recovery}\label{S5.2.5}

\begin{figure}[H]
 \centering
\includegraphics[width=1\textwidth]{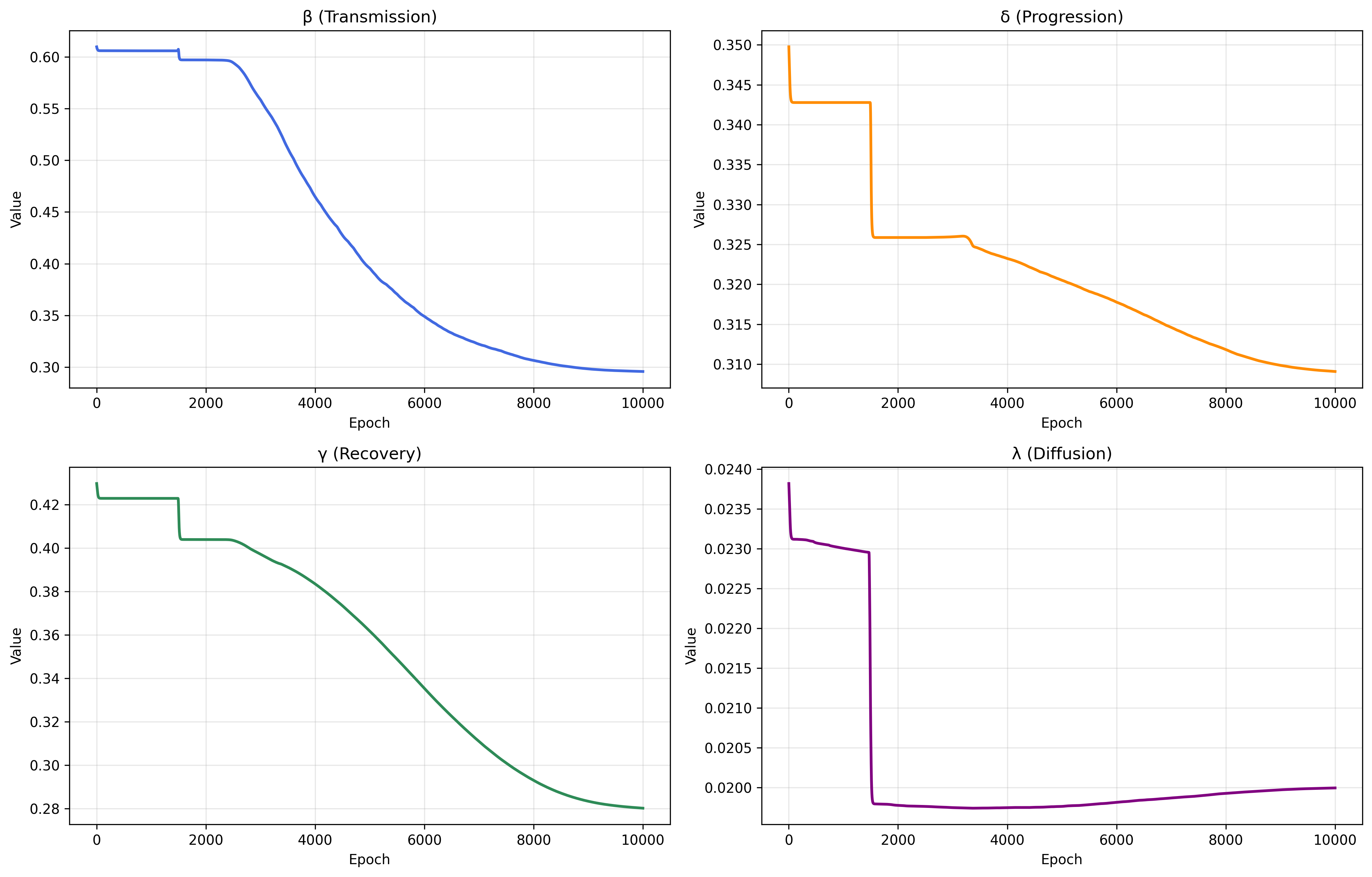}
 \caption{Convergence of estimated parameters during PINN training; 2D case.}\label{F11}
 \end{figure}

Figure~\ref{F11} complements this analysis by illustrating the convergence history of the four estimated parameters $\beta$, $\delta$, $\gamma$, and $\lambda$ in the 2D setting. 
The parameter recovery is slightly less accurate than in 1D case, as expected due to the larger hypothesis space and the smaller relative contribution of each data point to the parameter gradients. 
Table~\ref{Tab8} reports the comparison between the true and PINN-estimated values of the inferred parameters.

 \clearpage



\subsubsection{Quantitative Error Analysis}\label{S5.2.6}

Tables~\ref{Tab7} and~\ref{Tab8} summarize the pointwise error metrics for each compartment and the parameter recovery results, respectively.

{
\begin{table}[H]
\centering
\setlength{\tabcolsep}{0.7cm}
\caption{Quantitative error metrics between the PINN predictions and the reference data; 2D case.}\label{Tab7}
\adjustbox{max width=\textwidth}{
\begin{tabular}{ccccc}
\hline
\textbf{Compartment} & \textbf{Rel.\ $L^2$ Error} &
\textbf{MAE} & \textbf{RMSE} & \textbf{Max Error} \\
\hline\hline
$S$ & 3.1121\,e-03 & 2.3606\,e-03 & 3.0792\,e-03 & 9.8181\,e-03 \\
$E$ & 8.3622\,e-02 & 2.7817\,e-04 & 4.3298\,e-04 & 3.7180\,e-03 \\
$I$ & 2.8472\,e-01 & 1.3021\,e-03 & 1.7549\,e-03 & 4.6044\,e-03 \\
$R$ & 5.2685\,e-01 & 1.2035\,e-03 & 1.5878\,e-03 & 3.3886\,e-03 \\
\hline
\end{tabular}
}
\end{table}
\begin{table}[H]
\centering
\setlength{\tabcolsep}{0.7cm}
\caption{Comparison between the true and PINN-estimated values of the inferred parameters; 2D case.}\label{Tab8}
\adjustbox{max width=\textwidth}{
\begin{tabular}{ccccc}
\hline
\textbf{Parameter} & \textbf{True Value} & \textbf{Estimated} &
\textbf{Abs.\ Error} & \textbf{Rel.\ Error} \\
\hline\hline
$\beta$   & 0.40 & 0.2957 & 1.0426\,e-01 & 26.06\% \\
$\delta$  & 0.30 & 0.4091 & 1.0910\,e-01 & 36.37\% \\
$\gamma$  & 0.20 & 0.2802 & 8.0200\,e-02 & 40.10\% \\
$\lambda$ & 0.05 & 0.0200 & 3.0000\,e-02 & 60.00\% \\
\hline
\end{tabular}
}
\end{table}
}


In this case, the observed estimation errors are primarily due to the initialization and optimization of the unknown parameters within the broad interval  $(0,1)$ during the training process, 
which corresponds to a general identification setting.
If prior information about narrower admissible parameter intervals were available, the accuracy of the parameter estimation could be significantly improved. 
In the present work, however, we intentionally considered the general interval $(0,1)$ since no precise prior bounds on the parameters were assumed.

\section{Conclusion and Future Work}\label{S6}
In this work, we developed a constraint-aware PINN framework for solving and identifying parameters in SEIR reaction-diffusion epidemic systems with vital dynamics and homogeneous Neumann boundary conditions. 
The proposed methodology combines epidemiological constraints, observational data, and governing partial differential equations within a unified optimization framework capable of addressing both forward simulation and inverse parameter estimation problems.
To overcome the scarcity of high-resolution spatial epidemiological datasets, synthetic benchmark data were generated using positivity-preserving NSFD schemes. 
These schemes preserve essential qualitative properties of the epidemic dynamics, including non-negativity, boundedness, and numerical stability, thereby providing reliable reference solutions for training and validation.
The numerical experiments conducted in one- and two-dimensional spatial domains demonstrated that the proposed PINN framework accurately reconstructs the spatiotemporal evolution of all epidemiological compartments while maintaining low approximation errors throughout the computational domain. 
The proposed framework establishes an effective connection between mechanistic epidemiological modeling and modern scientific machine learning, offering promising perspectives for the analysis and prediction of complex spatiotemporal disease dynamics.


Several research directions may further extend the present work. 
First, the proposed framework may be generalized to more complex epidemiological structures, including SEIARD, age-structured, delayed, fractional-order, stochastic, or multi-strain reaction–diffusion systems. 
Incorporating heterogeneous diffusion coefficients and spatially varying transmission parameters would further improve the realism of the model in practical applications.
Another promising direction concerns the integration of real epidemiological datasets. 
Future studies may investigate the assimilation of noisy and incomplete observations obtained from mobility networks, geographic surveillance systems, or public health databases. 
Such developments would allow the framework to be applied to realistic forecasting, monitoring, and epidemic control scenarios.
In addition, extending the framework to stochastic reaction–diffusion epidemic models would enable the analysis of environmental variability, demographic fluctuations, random mobility patterns, and uncertainty propagation in spatiotemporal disease transmission.
From a computational perspective, future work may focus on improving the efficiency and scalability of PINN training for large-scale multidimensional problems. 
Possible directions include adaptive residual sampling strategies, domain decomposition techniques, operator-learning architectures, and hybrid numerical–neural methodologies.
Finally, incorporating optimal control strategies and uncertainty quantification within the PINN framework represents another important research direction. 
Such developments could provide valuable tools for evaluating intervention policies, estimating confidence intervals for inferred parameters, and supporting data-driven decision-making in epidemic management and public health planning.



\section*{Declarations}

\subsection*{Author Contributions}
\noindent
A. Zinihi: Conceptualization, Soft\-ware, Methodology, Validation, Formal Analysis, Investigation, Writing-Original Draft, Writing-Review and Editing, Visualization.\\
M. Ehrhardt: Conceptualization, Supervision, Methodology, Formal Analysis, Investigation, Writing-Review and Editing.

\subsection*{Computational Setup} 
\noindent
All numerical simulations for the PINN framework were conducted on a workstation running Microsoft Windows 11 Home. The system features an Intel Core Ultra 9 275HX CPU (24 cores, 24 threads) with 32 GB of RAM.\\
Computations were accelerated using a hybrid GPU setup, including an integrated Intel GPU and an NVIDIA GeForce RTX 5070 Ti Laptop GPU (12 GB VRAM). 

\subsection*{Software Environment} 
\noindent
The simulations were implemented in \textit{Python 3.13} using PyTorch.

\subsection*{Data Availability} 
\noindent
This study uses synthetic datasets as described in Sections~\ref{S3} and~\ref{S5}.
All data analyzed or generated in this work, which support the results, are included in this article.

\subsection*{Conflict of Interest} 
\noindent
The authors declare that there are no problems or conflicts of interest between them that may affect the study in this paper.


\bibliographystyle{unsrt}
\bibliography{References}

@article{Kermack1927,
 author = {Kermack, W. O. and McKendrick, A. G.},
 title = {A contribution to the mathematical theory of epidemics},
 fjournal = {Proceedings of the Royal Society A},
 journal = {Proc. Royal Soc. A},
doi = {10.1098/rspa.1927.0118},
 url = {https://doi.org/10.1098/rspa.1927.0118},
 year = {1927},
 publisher = {The Royal Society},
 volume = {115},
 number = {772},
 pages = {700--721}
}

@article{Hethcote2000,
  title = {The mathematics of infectious diseases},
  volume = {42},
  ISSN = {1095-7200},
  url = {http://dx.doi.org/10.1137/S0036144500371907},
  DOI = {10.1137/s0036144500371907},
  number = {4},
  journal = {SIAM Review},
  publisher = {Society for Industrial & Applied Mathematics (SIAM)},
  author = {Hethcote,  H. W.},
  year = {2000},
  pages = {599–653}
}

@article{Zinihi2025FDE,
  author = {Zinihi, A. and Sidi Ammi, M. R. and Torres, D. F. M.},
  title = {Fractional differential equations of a reaction-diffusion {SIR} model involving the {C}aputo-fractional time-derivative and a nonlinear diffusion operator},
  volume = {14},
  ISSN = {2163-2480},
  url = {http://dx.doi.org/10.3934/eect.2025018},
  DOI = {10.3934/eect.2025018},
  number = {5},
  fjournal = {Evolution Equations and Control Theory},
  journal = {Evolut. Eqs. Contr. Theo.},
  publisher = {American Institute of Mathematical Sciences (AIMS)},
  year = {2025},
  pages = {944–967}
}

@article{Wang2021,
author = {Wang,  N. and Zhang,  L. and Teng,  Z.},
  title = {Dynamics in a reaction-diffusion epidemic model via environmental driven infection in heterogenous space},
  volume = {16},
  ISSN = {1751-3766},
  url = {http://dx.doi.org/10.1080/17513758.2021.1900428},
  DOI = {10.1080/17513758.2021.1900428},
  number = {1},
  fjournal = {Journal of Biological Dynamics},
  journal = {J. Biol. Dynam.},
  publisher = {Informa UK Limited},
  year = {2021},
  pages = {373–396}
}

@article{Zinihi2025OC,
   author = {Zinihi, A. and Sidi Ammi, M. R. and Ehrhardt, M.},
  title = {Optimal control of a diffusive epidemiological model involving the {C}aputo-{F}abrizio fractional time-derivative},
  volume = {14},
  ISSN = {2666-8181},
  url = {http://dx.doi.org/10.1016/j.padiff.2025.101188},
  DOI = {10.1016/j.padiff.2025.101188},
  fjournal = {Partial Differential Equations in Applied Mathematics},
  journal = {Part. Diff. Eqs. Appl. Math.},
  publisher = {Elsevier BV},
  year = {2025},
  pages = {101188}
}

@book{Strikwerda2004,
  title = {Finite difference schemes and partial differential equations},
  ISBN = {9780898717938},
  url = {http://dx.doi.org/10.1137/1.9780898717938},
  DOI = {10.1137/1.9780898717938},
  publisher = {Society for Industrial and Applied Mathematics},
  author = {Strikwerda, J. C.},
  year = {2004}
}

@book{Banks1989,
  title={Estimation techniques for distributed parameter systems},
  ISBN = {9781461237006},
  url = {http://dx.doi.org/10.1007/978-1-4612-3700-6},
  DOI = {10.1007/978-1-4612-3700-6},
  publisher = {Birkh\"{a}user Boston},
  author = {Banks, H. T. and Kunisch, K.},
  year = {1989}
}

@article{Raissi2019,
  title = {Physics-informed neural networks: {A} deep learning framework for solving forward and inverse problems involving nonlinear partial differential equations},
  volume = {378},
  ISSN = {0021-9991},
  url = {http://dx.doi.org/10.1016/j.jcp.2018.10.045},
  DOI = {10.1016/j.jcp.2018.10.045},
  fjournal = {Journal of Computational Physics},
  journal = {J. Comput. Phys.},
  publisher = {Elsevier BV},
  author = {Raissi, M. and Perdikaris, P. and Karniadakis, G. E.},
  year = {2019},
  pages = {686–707}
}

@article{Karniadakis2021,
  author = {Karniadakis, G. E. and Kevrekidis, I. G. and Lu, L. and Perdikaris, P. and Wang, S. and Yang, L.},
  title = {Physics-informed machine learning},
  volume = {3},
  ISSN = {2522-5820},
  url = {http://dx.doi.org/10.1038/s42254-021-00314-5},
  DOI = {10.1038/s42254-021-00314-5},
  number = {6},
  fjournal = {Nature Reviews Physics},
  journal = {Nature Rev. Phys.},
  publisher = {Springer Science and Business Media LLC},
  year = {2021},
  pages = {422-440}
}

@article{Kissas2020,
author = {Kissas, G. and Yang, Y. and Hwuang, E. and Witschey, W. R. and Detre, J. A. and Perdikaris, P.},
  title = {Machine learning in cardiovascular flows modeling: {P}redicting arterial blood pressure from non-invasive {4D} flow {MRI} data using physics-informed neural networks},
  volume = {358},
  ISSN = {0045-7825},
  url = {http://dx.doi.org/10.1016/j.cma.2019.112623},
  DOI = {10.1016/j.cma.2019.112623},
  fjournal = {Computer Methods in Applied Mechanics and Engineering},
  journal = {Comput. Meth. Appl. Mech. Engrg.},
  publisher = {Elsevier BV},
  year = {2020},
  pages = {112623}
}

@article{Cuomo2022,
  author = {Cuomo, S. and Di Cola, V. S. and Giampaolo, F. and Rozza, G. and Raissi, M. and Piccialli, F.},
  title = {Scientific machine learning through physics–informed neural networks: {W}here we are and what’s next},
  volume = {92},
  ISSN = {1573-7691},
  url = {http://dx.doi.org/10.1007/S10915-022-01939-Z},
  DOI = {10.1007/s10915-022-01939-z},
  number = {3},
  fjournal = {Journal of Scientific Computing},
  journal = {J. Sci. Comput.},
  publisher = {Springer Science and Business Media LLC},
  year = {2022}
}

@article{Toscano2025,
author={Toscano, J. D. and Oommen, V. and Varghese, A. J. and Zou, Z. and Ahmadi Daryakenari, N. and Wu, C. and Karniadakis, G. E.},
  title={From {PINNs} to {PIKANs}: {R}ecent advances in physics-informed machine learning},
  volume = {1},
  ISSN = {3005-1436},
  url = {http://dx.doi.org/10.1007/s44379-025-00015-1},
  DOI = {10.1007/s44379-025-00015-1},
  number = {1},
  fjournal = {Machine Learning for Computational Science and Engineering},
  journal = {Machine Learn. Comput. Sci. Engrg.},
  publisher = {Springer Science and Business Media LLC},
  year = {2025}
}

@article{Zinihi2025MC,
author = {Zinihi, A. and Sidi Ammi, M. R. and Bachir, A.},
  title = {Multi-city modeling of epidemics using a topology-based {SIR} model: {Neural} network-enhanced {SAIRD} model},
  volume = {92},
  ISSN = {1877-7503},
  url = {http://dx.doi.org/10.1016/j.jocs.2025.102721},
  DOI = {10.1016/j.jocs.2025.102721},
  fjournal = {Journal of Computational Science},
  journal = {J. Comput. Sci.},
  publisher = {Elsevier BV},
  year = {2025},
  pages = {102721}
}

@article{Berkhahn2022,
author={Berkhahn, S. and Ehrhardt, M.},
  title={A physics-informed neural network to model {COVID-19} infection and hospitalization scenarios},
  volume = {2022},
  ISSN = {2731-4235},
  url = {http://dx.doi.org/10.1186/s13662-022-03733-5},
  DOI = {10.1186/s13662-022-03733-5},
  number = {1},
  fjournal = {Advances in Continuous and Discrete Models},
  journal = {Adv. Contin. Discr. Mod.},
  publisher = {Springer Science and Business Media LLC},
  year = {2022}
}

@article{Millevoi2024,
 author={Millevoi, C. and Pasetto, D. and Ferronato, M.},
  title={A physics-informed neural network approach for compartmental epidemiological models},
  volume = {20},
  ISSN = {1553-7358},
  url = {http://dx.doi.org/10.1371/journal.pcbi.1012387},
  DOI = {10.1371/journal.pcbi.1012387},
  number = {9},
  fjournal = {PLOS Computational Biology},
  journal = {PLOS Comput. Biol.},
  publisher = {Public Library of Science (PLoS)},
  editor = {Scarpino,  Samuel V.},
  year = {2024},
  pages = {e1012387}
}

@article{Nelson2024,
  title={Modeling the dynamics of {COVID-19} in {Japan}: {E}mploying data-driven deep learning approach},
  ISSN = {1868-808X},
  url = {http://dx.doi.org/10.1007/s13042-024-02301-5},
  DOI = {10.1007/s13042-024-02301-5},
  fjournal = {International Journal of Machine Learning and Cybernetics},
   journal = {Int. J. Mach. Learn.  Cybern.},
  publisher = {Springer Science and Business Media LLC},
  author={Nelson, S. P. and Raja, R.  and Eswaran, P. and Alzabut, J. and Rajchakit, G.},
  year = {2024}
}

@article{Shamsara2025,
  title={An informed deep learning model of the {Omicron} wave and the impact of vaccination},
  volume = {191},
  ISSN = {0010-4825},
  url = {http://dx.doi.org/10.1016/j.compbiomed.2025.109968},
  DOI = {10.1016/j.compbiomed.2025.109968},
  fjournal = {Computers in Biology and Medicine},
  journal = {Comput. Biol. Med.},
  publisher = {Elsevier BV},
  author={Shamsara, E. and K{\"o}nig, F. and Pfeifer, N.},
  year = {2025},
  pages = {109968}
}

@misc{Zinihi2025PINN, 
 doi = {10.48550/ARXIV.2509.22760}, 
 url = {https://arxiv.org/abs/2509.22760}, 
 author = {Zinihi, A.}, 
 title = {Identifying memory effects in epidemics via a fractional {SEIRD} model and physics-informed neural networks}, 
 note = {arXiv preprint}, 
 year = {2025}, 
 copyright = {Creative Commons Attribution 4.0 International}
}

@article{Igiri2025,
  title = {Monkeypox transmission dynamics using fractional disease informed neural network: {A} global and continental analysis},
  volume = {13},
  ISSN = {2169-3536},
  url = {http://dx.doi.org/10.1109/ACCESS.2025.3559005},
  DOI = {10.1109/access.2025.3559005},
  journal = {IEEE Access},
  publisher = {Institute of Electrical and Electronics Engineers (IEEE)},
  author={Igiri, C. P. and Shikaa, S.},
  year = {2025},
  pages = {77611–77641}
}

@article{Heldmann2023,
  author={Heldmann, F. and Berkhahn, S. and Ehrhardt, M. and Klamroth, K.},
  title={{PINN} training using biobjective optimization: {T}he trade-off between data loss and residual loss},
  volume = {488},
  ISSN = {0021-9991},
  url = {http://dx.doi.org/10.1016/j.jcp.2023.112211},
  DOI = {10.1016/j.jcp.2023.112211},
  fjournal = {Journal of Computational Physics},
  journal = {J. Comput. Phys.},
  publisher = {Elsevier BV},
  year = {2023},
  pages = {112211}
}

@article{Cheng2024,
author={Cheng, H. and Mao, Y. and Jia, X.},
  title={A framework based on physics-informed graph neural {ODE}: {F}or continuous spatial-temporal pandemic prediction},
  volume = {54},
  ISSN = {1573-7497},
  url = {http://dx.doi.org/10.1007/s10489-024-05834-y},
  DOI = {10.1007/s10489-024-05834-y},
  number = {24},
  fjournal = {Applied Intelligence},
  journal = {Appl. Intell.},
  publisher = {Springer Science and Business Media LLC},
  year = {2024},
  pages = {12661–12675}
}

@article{Liu2024,
  title={A combined neural {ODE-Bayesian} optimization approach to resolve dynamics and estimate parameters for a modified {SIR} model with immune memory},
  volume = {10},
  ISSN = {2405-8440},
  url = {http://dx.doi.org/10.1016/j.heliyon.2024.e38276},
  DOI = {10.1016/j.heliyon.2024.e38276},
  number = {19},
  journal = {Heliyon},
  publisher = {Elsevier BV},
  author={Liu, D. and Sopasakis, A.},
  year = {2024},
  pages = {e38276}
}

@article{Shan2025,
  title={{VI-PINNs}: {V}ariance-involved physics-informed neural networks for fast and accurate prediction of partial differential equations},
  volume = {623},
  ISSN = {0925-2312},
  url = {http://dx.doi.org/10.1016/j.neucom.2025.129360},
  DOI = {10.1016/j.neucom.2025.129360},
  journal = {Neurocomputing},
  publisher = {Elsevier BV},
  author={Shan, B. and Li, Y. and Huang, S.-J.},
  year = {2025},
  pages = {129360}
}

@article{Madden2024,
 author = {Madden, W. G. and Jin, W. and Lopman, B. and Zufle, A. 
 and Dalziel, B. and Metcalf,  C. Jessica E.  and Grenfell, B. T. and Lau, M. S. Y.},
  title = {Deep neural networks for endemic measles dynamics: {C}omparative analysis and integration with mechanistic models},
  volume = {20},
  ISSN = {1553-7358},
  url = {http://dx.doi.org/10.1371/journal.pcbi.1012616},
  DOI = {10.1371/journal.pcbi.1012616},
  number = {11},
  fjournal = {PLOS Computational Biology},
  journal = {PLOS Comput. Biol.},
  publisher = {Public Library of Science (PLoS)},
  year = {2024},
  pages = {e1012616}
}

@misc{Zinihi2025S,
  doi = {10.48550/ARXIV.2507.09328},
  url = {https://arxiv.org/abs/2507.09328},
  author = {Zinihi, A. and Ehrhardt, M. and Sidi Ammi, M. R.},
  title = {Spatiotemporal {SEIQR} epidemic modeling with optimal control for vaccination, treatment, and social measures},
  note = {arXiv preprint},
  year = {2025},
  copyright = {Creative Commons Attribution 4.0 International}
}

@article{Chang2022,
 author = {Chang, L. and Gao, S. and Wang, Z.},
 title = {Optimal control of pattern formations for an {SIR} reaction-diffusion epidemic model},
  volume = {536},
  ISSN = {0022-5193},
  url = {http://dx.doi.org/10.1016/j.jtbi.2022.111003},
  DOI = {10.1016/j.jtbi.2022.111003},
  fjournal = {Journal of Theoretical Biology},
   journal = {J. Theor. Bio.},
   publisher = {Elsevier BV},
   year = {2022},
  pages = {111003}
}

@article{Wang2012,
  author = {Wang, W. and Cai, Y. and Wu, M. and Wang, K. and Li, Z.},
  title = {Complex dynamics of a reaction-diffusion epidemic model},
  volume = {13},
  ISSN = {1468-1218},
  url = {http://dx.doi.org/10.1016/j.nonrwa.2012.01.018},
  DOI = {10.1016/j.nonrwa.2012.01.018},
  number = {5},
  fjournal = {Nonlinear Analysis: Real World Applications},
  journal = {Nonlin. Anal. Real World Appl.},
  publisher = {Elsevier BV},
  year = {2012},
  pages = {2240-2258}
}

@article{Chen2025,
  title = {A deep learning framework for predicting the spread of diffusion diseases},
  volume = {33},
  ISSN = {2688-1594},
  url = {http://dx.doi.org/10.3934/era.2025110},
  DOI = {10.3934/era.2025110},
  number = {4},
  journal = {Electronic Research Archive},
  publisher = {American Institute of Mathematical Sciences (AIMS)},
  author = {Chen, X. and Li, F. and Lian, H. and Wang, P.},
  year = {2025},
  pages = {2475-2502}
}

@article{Ye2025,
  title = {Integrating artificial intelligence with mechanistic epidemiological modeling: {A} scoping review of opportunities and challenges},
  volume = {16},
  ISSN = {2041-1723},
  url = {http://dx.doi.org/10.1038/s41467-024-55461-x},
  DOI = {10.1038/s41467-024-55461-x},
  number = {1},
  journal = {Nature Communications},
  publisher = {Springer Science and Business Media LLC},
  author={Ye, Y. and Pandey, A. and Bawden, C. and Sumsuzzman, D. M. and Rajput, R. and Shoukat, A. and Singer, B. H. and Moghadas, S. M. and Galvani, A. P.},
  year = {2025}
}

@book{Mickens1993,
  title = {Nonstandard finite difference models of differential equations},
  ISBN = {9789814440882},
  url = {http://dx.doi.org/10.1142/2081},
  DOI = {10.1142/2081},
  publisher = {World Scientific},
  author = {Mickens, R. E.},
  year = {1993}
}

@book{Mickens2020,
  title={Nonstandard finite difference schemes: {M}ethodology and applications},
  author={Mickens, R. E.},
  isbn={9789811222559},
  url={https://books.google.co.ma/books?id=PZwPEAAAQBAJ},
  year={2020},
  publisher={World Scientific Publishing Company}
}

@article{Zinihi2025NSFD,
author = {Zinihi, A. and Ehrhardt, M. and Sidi Ammi, M. R.},
  title = {A nonstandard finite difference scheme for an {SEIQR} epidemiological {PDE} model},
  volume = {520},
  ISSN = {0096-3003},
  url = {http://dx.doi.org/10.1016/j.amc.2026.129953},
  DOI = {10.1016/j.amc.2026.129953},
  fjournal = {Applied Mathematics and Computation},
  journal = {Appl. Math. Comput.},
  publisher = {Elsevier BV},
  year = {2026},
  pages = {129953}
}

@article{Ehrhardt2013,
  title = {A nonstandard finite difference scheme for convection-diffusion equations having constant coefficients},
  volume = {219},
  ISSN = {0096-3003},
  url = {http://dx.doi.org/10.1016/j.amc.2012.12.068},
  DOI = {10.1016/j.amc.2012.12.068},
  number = {12},
  fjournal = {Applied Mathematics and Computation},
  journal = {Appl. Math. Comput.},
  publisher = {Elsevier BV},
  author = {Ehrhardt, M. and Mickens, R. E.},
  year = {2013},
  pages = {6591-6604}
}

@article{ChenCharpentier2013,
author = {Chen-Charpentier, B. M. and Kojouharov, H. V.},
  title = {An unconditionally positivity preserving scheme for advection-diffusion reaction equations},
  volume = {57},
  ISSN = {0895-7177},
  url = {http://dx.doi.org/10.1016/j.mcm.2011.05.005},
  DOI = {10.1016/j.mcm.2011.05.005},
  number = {9–10},
  fjournal = {Mathematical and Computer Modelling},
  journal = {Math. Comput. Model.},
  publisher = {Elsevier BV},
  year = {2013},
  pages = {2177–2185}
}

@article{Pasha2023,
  title = {On the nonstandard finite difference method for reaction–diffusion models},
  volume = {166},
  ISSN = {0960-0779},
  url = {http://dx.doi.org/10.1016/j.chaos.2022.112929},
  DOI = {10.1016/j.chaos.2022.112929},
  journal = {Chaos,  Solitons \& Fractals},
  publisher = {Elsevier BV},
  author = {Pasha, S. A. and Nawaz, Y. and Arif, M. S.},
  year = {2023},
  pages = {112929}
}

@article{Wang2021S,
  title = {Understanding and mitigating gradient pathologies in physics-informed neural networks},
  volume = {43},
  ISSN = {1095-7197},
  url = {http://dx.doi.org/10.1137/20M1318043},
  DOI = {10.1137/20m1318043},
  number = {5},
  fjournal = {SIAM Journal on Scientific Computing},
  journal = {SIAM J. Sci. Comput.},
  publisher = {Society for Industrial \& Applied Mathematics (SIAM)},
  author = {Wang, S. and Teng, Y. and Perdikaris, P.},
  year = {2021},
  pages = {A3055–A3081}
}

@inproceedings{Tancik2020,
  title={Fourier features let networks learn high frequency functions in low dimensional domains},
  author={Tancik, M. and Srinivasan, P. P. and Mildenhall, B. and Fridovich-Keil, S. and Raghavan, N. and Singhal, U. and Ramamoorthi, R. and Barron, J. T. and Ng, R.},
  booktitle={Advances in Neural Information Processing Systems (NeurIPS)},
  volume={33},
  pages={7537--7547},
  year={2020}
}

@inproceedings{Glorot2010,
  title={Understanding the difficulty of training deep feedforward neural networks},
  author={Glorot, X. and Bengio, Y.},
  booktitle={Proceedings of the 13th international conference on artificial intelligence and statistics},
  pages={249-256},
  year={2010},
  organization={JMLR Workshop and Conference Proceedings}
}

@article{Liu1989,
  title = {On the limited memory {BFGS} method for large scale optimization},
  volume = {45},
  ISSN = {1436-4646},
  url = {http://dx.doi.org/10.1007/BF01589116},
  DOI = {10.1007/bf01589116},
  number = {1-3},
  fjournal = {Mathematical Programming},
  journal = {Math. Programm.},
  publisher = {Springer Science and Business Media LLC},
  author = {Liu, D. C. and Nocedal, J.},
  year = {1989},
  pages = {503–528}
}

@misc{UKGov2025,
  author = {UK Health Security Agency},
  title = {{Epidemiology of COVID-19 in England: January 2020 to December 2024}; \url{https://www.gov.uk/government/publications/epidemiology-of-covid-19-in-england/epidemiology-of-covid-19-in-england-january-2020-to-december-2024}},
  year = {Accessed on: May 10, 2026},
  url = {https://www.gov.uk/government/publications/epidemiology-of-covid-19-in-england/epidemiology-of-covid-19-in-england-january-2020-to-december-2024}
}

@misc{WHO2024Covid,
  author = {World Health Organization},
  title = {{COVID}-19 epidemiological update, 24 {December} 2024; \url{https://www.who.int/publications/m/item/covid-19-epidemiological-update---24-december-2024}},
  year = {Accessed on: May 10, 2026},
  url = {https://www.who.int/publications/m/item/covid-19-epidemiological-update---24-december-2024}
}

@misc{CDCFlu,
  author = {Centers for Disease Control and Prevention},
  title = {Seasonal Influenza ({Flu}), \url{https://www.cdc.gov/flu/index.htm}}, 
  year = {Accessed on: May 11, 2026},
  url = {https://www.cdc.gov/flu/index.htm}
}

@misc{ECDCFlu,
  author = {European Centre for Disease Prevention and Control},
  title = {Seasonal influenza data; \url{https://www.ecdc.europa.eu/en/seasonal-influenza}},
  year = {Accessed on: May 11, 2026},
  url = {https://www.ecdc.europa.eu/en/seasonal-influenza}
}

@misc{JHU2020,
  author = {Johns Hopkins University},
  title = {{COVID-19} dashboard; \url{https://coronavirus.jhu.edu/data}},
  year = {Accessed on: May 11, 2026 (2020)},
  url = {\url{https://coronavirus.jhu.edu/data}}
}

@article{Zinihi2026Koopman,
  title = {A {Koopman} operator framework for nonlinear epidemic dynamics: {Application} to an {SIRSD} model},
  volume = {6},
  ISSN = {3050-5178},
  url = {http://dx.doi.org/10.1016/j.nls.2025.100095},
  DOI = {10.1016/j.nls.2025.100095},
  journal = {Nonlinear Science},
  publisher = {Elsevier BV},
  author = {Zinihi, A. and Ehrhardt, M. and Sidi Ammi, M. R.},
  year = {2026},
  month = Apr,
  pages = {100095}
}

@misc{CDCCOVID,
  author = {Centers for Disease Control and Prevention},
  title = {{CDC COVID} data tracker; \url{https://www.cdc.gov/covid-data-tracker/index.html}},
  year = {Accessed on: May 11, 2026},
  url = {\url{https://www.cdc.gov/covid-data-tracker/index.html}}
}

@article{Liu2022,
  author = {Liu, P. and Li, H. X.},
  title = {Global behavior of a multi-group {SEIR} epidemic model with spatial diffusion in a heterogeneous environment},
  volume = {32},
  ISSN = {2083-8492},
  url = {http://dx.doi.org/10.34768/amcs-2022-0020},
  DOI = {10.34768/amcs-2022-0020},
  number = {2},
  journal = {International Journal of Applied Mathematics and Computer Science},
  publisher = {University of Zielona Góra, Poland},
  year = {2022}
}

@book{pazy2012semigroups,
  author={Pazy, A.},
  title={Semigroups of linear operators and applications to partial differential equations},
  year={2012},
  publisher={Springer Science \& Business Media}
}

@article{maamar2024nonstandard,
  title={A nonstandard finite difference scheme for a time-fractional model of {Zika} virus transmission},
  author={Maamar, M. H. and Ehrhardt, M. and Tabharit, L.},
  fjournal={Mathematical Biosciences and Engineering},
  journal={Math. Biosci. Engrg.},
  volume={21},
  number={1},
  pages={924--962},
  year={2024},
  publisher={American Institute of Mathematical Sciences}
}

@article{wang2022,
  title={Is {$L^2$} Physics Informed Loss Always Suitable for Training Physics Informed Neural Network?},
  author={Wang, C. and Li, S. and He, D. and Wang, L.},
  fjournal={Advances in Neural Information Processing Systems},
  journal={Adv. Neural Inform. Process. Syst.},
  volume={35},
  pages={8278--8290},
  year={2022}
}


\end{document}